
\input  amstex
\input amsppt.sty
\magnification1200
\vsize=23.5truecm
\hsize=16.5truecm
\NoBlackBoxes

\def\rn{{\Bbb R}^n}
\def\rnp{{\Bbb R}^n_+}
\def\rnm{\Bbb R^n_-}

\def\crnp{\overline{\Bbb R}^n_+}
\def\crnm{\overline{\Bbb R}^n_-}

\def\comega{\overline\Omega }
\def\ang#1{\langle {#1} \rangle}
\def\rp{ \Bbb R_+}

\def\Op{\operatorname{Op}}
\def\N{\Bbb N}
\def\R{\Bbb R}
\def\C{\Bbb C}
\def\Z{\Bbb Z}
\def\ol{\overline}
\def\E{\Cal E}
\def\F{\Cal F}
\def\simto{\overset\sim\to\rightarrow}
\def\supp{\operatorname{supp}}

\document
\topmatter
\title
The principal transmission condition 
\endtitle
\author Gerd Grubb \endauthor
\affil
{Department of Mathematical Sciences, Copenhagen University,\\
Universitetsparken 5, DK-2100 Copenhagen, Denmark.\\
E-mail {\tt grubb\@math.ku.dk}}\endaffil
\rightheadtext{Principal transmission}

\abstract
The paper treats  pseudodifferential operators $P=\Op(p(\xi ))$ with
homogeneous complex symbol $p(\xi )$ of
order $2a>0$, generalizing  the fractional Laplacian $(-\Delta
)^a$ but  lacking its symmetries, and taken to act on the halfspace $\rnp$.
The operators are seen to satisfy a
 principal $\mu $-transmission condition
relative to $\rnp$, but generally not the full $\mu $-transmission
condition satisfied by $(-\Delta )^a$ and related operators (with $\mu
=a$).
However, $P$
acts well on the so-called $\mu $-transmission
spaces over $\rnp$ (defined in earlier works), and
when $P$ moreover is
strongly elliptic, these spaces are the solution spaces for the
homogeneous Dirichlet problem for $P$, leading to 
regularity results with a factor $x_n^\mu $
(in a limited range of Sobolev spaces).
The information is then shown
to be sufficient to establish an 
integration by parts formula over $\rnp$ for $P$ acting on such
functions. The formulation in Sobolev spaces, and the results on
strongly elliptic operators going beyond operators with real kernels,
are new. 

Furthermore, large solutions with nonzero Dirichlet traces are described, and a halfways Green's formula is
established, for this new class of operators.

Since the principal $\mu $-transmission condition has weaker
requirements than the
full $\mu $-transmission condition assumed in earlier papers, new
arguments were needed, relying on work of Vishik and Eskin instead of the
Boutet de Monvel theory. The results cover the case of
nonsymmetric operators with real kernel that
were only partially treated in a preceding paper.
\endabstract



\endtopmatter

\subhead 1. Introduction \endsubhead

Boundary value problems for fractional-order pseudodifferential
operators $P$, in particular where $P$ is a generalization of the
fractional Laplacian $(-\Delta )^a$ ($0<a<1$), have currently received much interest in
applications, such as in financial theory and probability (but
also in mathematical physics and differential geometry), and many
methods have been used, most often probabilistic or potential-theoretic methods.

The author has studied such problems by pseudodifferential methods in \cite{G15--G21}, under the assumption
that the operators satisfy a $\mu $-transmission condition at the boundary of
the domain $\Omega \subset \rn$, which allows to show regularity
results for solutions of  the Dirichlet problem in elliptic
cases, to show integration by parts formulas, and much else.

In the present paper we consider translation-invariant
pseudodifferential  operators \linebreak ($\psi $do's)
$P=\Op(p(\xi ))$ of order $2a>0$ with homogeneous symbol $p(\xi )$,
which are only taken to satisfy
the top-order equation in the $\mu $-transmission condition
(relative to the domain $\Omega =\rnp$), we
call this the {\it principal $\mu $-transmission condition}. It is shown
that they retain some of the features: The solution spaces for the
homogeneous Dirichlet problem in the elliptic case equal the $\mu
$-transmission spaces from \cite{G15} (in a setting of low-order Sobolev spaces), having a factor
$x_n^\mu $. The
integration by parts formula holds (even when $P$ is not
elliptic):
$$
\int_{\rnp}(  Pu\,\partial_n\bar
u'+\partial_nu\overline{ P^*u'})\,dx
=\Gamma (\mu +1){\Gamma(\mu '+1)}\int_{{\Bbb
R}^{n-1}}s_0\gamma _0(u/x_n^{\mu })\,{\gamma _0(\bar u'/x_n^{\mu '})}\,
dx',
$$
when $u$ and $u'$ are in $x_n^\mu C^\infty (\crnp)$ resp.\ $x_n^{\bar\mu
'}C^\infty (\crnp)$ ($\mu '=2a-\mu $) and compactly supported.
We also treat nonhomogeneous local Dirichlet problems with 
Dirichlet trace $\gamma _0(u/x_n^{\mu -1})$, and show how the
above formula implies a ``halfways'' Green's formula where one factor
has nonzero
Dirichlet trace. $P$ can be of any positive order, and $\mu $ can be complex.

The results apply in particular to the operator 
$L=\Op(\Cal A(\xi )+i\Cal B(\xi ))$ with $\Cal A$ real, positive and even
in $\xi $,  $\Cal B$ real and odd
in $\xi $,  which satisfies the principal $\mu $-transmission
equation for a suitable real $\mu $.
Hereby we can compensate for an error made in the recent publication
\cite{G21} (see also \cite{G22}), where it was overlooked that $L$ may not satisfy the full
$\mu $-transmission condition when $\Cal B\ne 0$ (it does so for $\Cal
B=0$). The general $L$
are now covered by the present work. They were treated earlier by
Dipierro, Ros-Oton, Serrra and Valdinoci \cite{DRSV21} under some
hypotheses on $a$ and $\mu $; they come up in applications as
infinitesimal
generators of  $\alpha $-stable $n$-dimensional L\'evy processes, see
\cite{DRSV21}. (The calculations in \cite{G21}
are valid when applied to operators satisfying the full $\mu
$-transmission condition.)

The study of $x$-independent $\psi $do's $P$ on the half-space $\rnp$
serves as a model case for $x$-dependent operators on domains $\Omega
\subset \rn$ with curved boundary, and can be expected to be a useful
ingredient in the general treatment, as carried out for the
operator $L$ in \cite{DRSV21}.

\medskip

\noindent{\it Plan of the paper:} In Section 2 we give an overview of
the aims and results of the paper with only few technicalities. Section 3 introduces the principal
transmission condition in detail for homogeneous $\psi $do symbols. In Section
4, the Wiener-Hopf method is applied to derive basic decomposition and
factorization formulas for such symbols. This is used in
Section 5 to establish mapping properties for the operators, and regularity properties for solutions of the
homogeneous Dirichlet problem in strongly elliptic cases; here $\mu
$-transmission spaces (known from \cite{G15}) defined in an $L_2$-framework play an important
role. Section 6 gives the proof of the above-mentioned integration by
parts formula on $\rnp$.
Section 7 treats nonhomogeneous local
Dirichlet conditions, and a halfways Green's formula is established.

\subhead 2. Presentation of the main results \endsubhead

The study is concerned with the so-called {\it model case}, where the
pseudodifferential operators have $x$-independent symbols, hence act as
simple multiplication operators in the Fourier transformed space (this frees
us from using the deeper composition rules needed for $x$-dependent
symbols), and the considered open subset $\Omega $ of $\rn$ is
simplest possible, namely $\Omega =\rnp = \{x\in \rn\mid x_n>0\}$. We
assume $n\ge 2$ and
denote $x=(x_1,\dots,x_n)=(x',x_n)$, $x'=(x_1,\dots,x_{n-1})$.
Recall the formulas for the Fourier transform $\F$ and the operator
$P=\Op(p(\xi ))$:
$$
\aligned
\Cal F
u&= \hat
u(\xi )= 
\int_{{\Bbb R}^n}e^{-ix\cdot \xi }u(x)\, dx,\quad \Cal F^{-1}v=(2\pi )^{-n}
\int_{{\Bbb R}^n}e^{ix\cdot \xi }v(\xi )\, d\xi ,\\
Pu&=\Op(p(\xi ))u=\Cal F^{-1}(p(\xi )(\Cal Fu)(\xi )).
\endaligned
\tag2.1
$$
We work
in $L_2(\rn)$ and $L_2(\rnp)$   and their derived $L_2$-Sobolev spaces
(the reader is urged to consult (5.1) below for notation). On $L_2(\rn)$, the Plancherel theorem
$$
\| u\|_{L_2(\rn)}=c\|\hat u\|_{L_2(\rn)},\quad c=(2\pi )^{-n/2},
\tag2.2
$$
makes norm estimates of operators easy. (There is more on  Fourier transforms
and distribution theory e.g.\ in \cite{G09}.) The model case serves both as a
simplified special case, and as a proof ingredient for more general
cases of  domains with curved boundaries, and possibly $x$-dependent
symbols.

The symbols $p(\xi )$ we shall consider are scalar and homogeneous of degree
$m=2a>0$ in $\xi $, i.e., $p(t\xi )=t^mp(\xi )$ for $t>0$, and are $C^1$ for $\xi \ne 0$, defining operators
$P=\Op(p)$. 

A typical example  is the squareroot Laplacian with drift:
$$
L_1=(-\Delta )^\frac12+\pmb b\cdot \nabla,\text{ with symbol }\Cal
L_1(\xi )=|\xi |+i\pmb b\cdot \xi ,\tag2.3
$$
where $\pmb b=(b_1,\dots, b_n)$ is a real vector. Here $m=1$,
$a=\frac12$. It satisfies the condition for  {\it strong
ellipticity}, which is:
$$
\operatorname{Re}p(\xi )\ge c_0|\xi |^m\text{ with }c_0>0;\tag2.4
$$
this is important in regularity discussions. Some results are obtained without the ellipticity hypothesis;
as an example we can take the operator $L_2$ with symbol
$$
\Cal L_2(\xi )=|\xi _1+\dots+\xi _n|+i\pmb b\cdot \xi ,\tag2.5
$$
whose real part is zero e.g.\ when $\xi=(1,-1,0,\dots,0)$.

The operators are well-defined on the Sobolev spaces over $\rn$: When
$p$ is homogeneous of degree $m\ge 0$, there is an inequality
$$
|p(\xi )|\le C|\xi |^m\le C\ang \xi ^m,\quad \ang \xi =(1+|\xi |^2)^\frac12
$$ 
(we say that $p$ {\it is of order} $m$); then 
$$
\|Pu\|_{L_2(\rn)}=c\|p(\xi )\hat u(\xi )\|_{L_2(\rn)}\le cC\|\ang\xi ^m\hat u\|_{L_2(\rn)}=C'\|u\|_{H^m(\rn)},\tag2.6
$$
so $P$ maps $H^m(\rn)$ continuously into $L_2(\rn)$. Similarly, it
maps $H^{s+m}(\rn)$ continuously into $H^s(\rn)$ for all $s\in\R$.

But for these pseudodifferential operators it is not obvious how to define them relative to
the subset $\rnp$, since they are not defined pointwise like
differential operators, but by
integrals (they are {\it nonlocal}). 
The convention is here to let them act on suitable linear subsets of $L_2(\rnp)$,
where we identify $L_2(\rnp)$ with the set  of $u\in L_2(
\rn)$ that are zero on $\rnm$, i.e., have their support
$\operatorname{supp}u\subset\crnp$. (The support $\operatorname{supp}u$ of a
function or distribution  $u$ is the complement of the largest open
set where $u=0$. The operator that extends
functions on $\rnp$ by zero
on $\rnm$ is denoted $e^+$.)
Then we apply $P$  and restrict
to $\rnp$ afterwards; this is the operator $r^+P$. ($r^+$ stands for
restriction to $\rnp$.)

Aiming for the integration by parts formula mentioned in the start, we
have to clarify for which functions $u,u'$ the integrals make sense.
It can be expected from earlier studies (\cite{RS14}, \cite{G16}, \cite{DRSV21})
that the integral will be meaningful for solutions of the 
so-called
{\it homogeneous Dirichlet  problem on} $\rnp$, namely the problem
$$
r^+Pu=f\text{ on }\rnp, \quad u=0 \text{ on }\rnm
\tag2.7
$$
(where the latter condition can also be written $ \supp u\subset  \crnp$).
This raises the question of where $r^+P$ lands; which $f$ can be
prescribed? Or, if $f$ is given in certain space, where should $u$ lie
in order to hit the
space where $f$ lies?

Altogether, we address the following three questions on $P$:
\roster
\item Forward mapping properties. From which spaces does $r^+P$ map 
into an $H^s$-space for $f$?
\item Regularity properties. If $u$ solves (2.7) with $f$ in an
$H^s$-space for a high $s$, will $u$ then belong to a space with a
similar high regularity?
\item Integration by parts formula for functions in spaces
where $r^+P$ is  well-defined.
\endroster

It turns out that the answers to all three points depend profoundly on the
introduction of so-called $\mu $-transmission spaces. To explain their
importance, we turn for a moment to the fractional Laplacian which has a
well-established treatment:

For the case of  $(-\Delta )^a$, $0<a<1$, it
was shown in \cite{G15} that the following space is relevant:
$$
\Cal E_a(\crnp)=e^+x_n^aC^\infty (\crnp).\tag 2.8
$$
It has the property that $(-\Delta )^a$ maps it to $C^\infty (\crnp)$;
more precisely,
$$
r^+(-\Delta )^a\text{ maps }\Cal E_a(\crnp)\cap \E'(\rn)\text{ into }C^\infty (\crnp).\tag2.9
$$
Here $\E'(\rn)$ is the space of distributions with compact support, so
the intersection with this space means that we consider functions in
$\E_a$ that are zero outside a compact set.

For Sobolev spaces, it was found in \cite{G15} that the good space for $u$ is the
so-called $a$-transmission space  $H^{a(t)}(\crnp)$; here
$$
r^+(-\Delta )^a\text{ maps }H^{a(t)}(\crnp)\cap \E'(\rn)\text{ into
}\ol H^{t-2a}(\rnp),\tag2.10
$$
for all $t\ge a$ (say). $\E_a(\crnp)\cap
\E'(\rn)$ is a dense subset of  $H^{a(t)}(\crnp)$.  The definition of the space $H^{a(t)}(\crnp)$ is recalled 
below in (2.15) and in more detail in Section 5.3; let us for the
moment  just mention that it is the sum of the space $\dot
H^t(\crnp)$ and a certain subspace of $x_n^a\ol H^{t-a}(\rnp)$. This
also holds when $a$ is replaced by a more general $\mu $.

For $(-\Delta )^a$, the $a$-transmission spaces provide the right answers
to question (1), and they are likewise right for question (2) (both
facts established in \cite{G15}), and there are integration by parts formulas for
$(-\Delta )^a$ applied to elements of these spaces, \cite{G16, G18}.

The key to the proofs is the so-called $a$-transmission condition that
$(-\Delta )^a$ satisfies; it is an infinite list of equations for
$p(\xi )$ and its derivatives, linking the values on the
interior normal to $\rnp$ with the values on the exterior normal. We
formulate it below with $a$ replaced
by a general $\mu $.

\proclaim{Definition 2.1} Let $\mu \in\C$, and let $p(\xi )$ be
homogeneous of degree $m$. Denote the interior resp.\ exterior
normal to the boundary of $\rnp$ by  $(0,\pm 1)=\{(\xi ',\xi
_n)\mid \xi '=0,\xi _n=\pm 1\}$.

$1^\circ$ $p$ (and $P=\Op(p)$) is said to  satisfy {\bf the principal $\mu $-transmission
condition} at $\rnp$ if
$$
p(0,-1)=e^{i\pi (m-2\mu )}p(0,1).\tag2.11
$$

$2^\circ$ $p$ (and $P=\Op(p)$) is said to  satisfy {\bf the  $\mu $-transmission
condition} at $\rnp$ if
$$
\partial_\xi ^\alpha p(0,-1)=e^{i\pi (m-2\mu -|\alpha |)}\partial_\xi
^\alpha p(0,1), \text{  for all }\alpha \in \N_0^n.\tag2.12
$$ 

\endproclaim

Note that $\mu $ is determined from $p$ in (2.11) up to addition of an
integer, when $p(0,1)\ne 0$.

The operators considered on smooth domains $\Omega $ in \cite{G15}
were assumed to satisfy (2.12) (for the top-order term $p_0$ in the symbol) at
all boundary points $x_0\in \partial\Omega $, with $(0,1)$ replaced by
the interior normal $\nu $ at $x_0$, and $(0,-1)$ replaced by $-\nu $. The lower-order terms $p_j$ in
the symbol, homogeneous of degree $m-j$, should then satisfy analogous
rules with $m-j$ instead of $m$.

The principal $\mu $-transmission condition (2.11) is of course much less 
demanding than the full $\mu $-transmission condition (2.12).
What we show in the present paper is
 that when (2.11) holds, the $\mu $-transmission
spaces are still relevant, and provide the appropriate answers to both
questions (1) and (2), however just for $t$ (the regularity parameter)
in a limited range. This range is large enough that integration by
parts formulas can be established, answering (3).

By simple geometric considerations one finds:

\proclaim{Proposition 2.2} $1^\circ$ When $p(\xi )$ is homogeneous of degree $m$,
there is a $\mu \in \C$, uniquely determined modulo $\Z$ if $p(0,1)\ne
0$, such that
{\rm (2.11)} holds.

$2^\circ$ If moreover, $p$ is strongly elliptic {\rm (2.4)} and $m=2a>0$, $\mu $
can be chosen uniquely to satisfy $\mu =a+\delta $ with $|\operatorname{Re}\delta |<\frac12$.
\endproclaim

This is shown in Section 3. From here on we work under  two slightly different
assumptions. The symbol  $p(\xi )$ is in both cases taken homogeneous of degree
$m=2a>0$ and $C^1$ for $\xi \ne 0$. We pose  Assumption 3.1 requiring that $p$
is strongly elliptic and $\mu $ is chosen as in Proposition 2.2 $2^\circ$.
We pose Assumption 3.2 just requiring that $\mu $ is defined according to
Proposition 2.2 $1^\circ$. In all cases we write $\mu =a+\delta $, and
define $\mu '=a-\delta =2a-\mu $.

\example{Example 2.3} Consider $\Cal L_1=|\xi |+i\pmb b\cdot \xi $
defined in (2.3). The order is $1=m=2a$, so $a=\frac12$.  Here $\Cal L_1(0,1)=1+ib_n$ and $\Cal L_1(0,-1)=1-ib_n$. The
angle $\theta $ in $\C=\R^2$ between the positive real axis and
$1+ib_n$ is $\theta =\operatorname{Arctan}b_n$. Set $\delta =\theta
/\pi $, then 
$$
\aligned
\Cal L_1(0,1)&=e^{i\pi \delta }|\Cal L_1(0,1)|=e^{i\pi \delta
}(1+|\pmb b|^2)^\frac12,
\text{ similarly }\\
\Cal L_1(0,-1)&=e^{-i\pi \delta }|\Cal L_1(0,-1)|=e^{-i\pi \delta
}(1+|{\pmb b}|^2)^\frac13.
\endaligned
$$
 Moreover, 
$$
\Cal L_1(0,-1)/\Cal L_1(0,1)=e^{-2i\pi \delta }=e^{i\pi (2a-2(a+\delta
))},\; a=\tfrac12,
$$
so (2.11) holds with $m=2a=1$, $\mu =\frac12+\delta $, where $\delta
=\frac1\pi \operatorname{Arctan}b_n$, and Assumption 3.1 is satisfied. Note that $\delta \in
\,]-\frac12,\frac12[\,$.

For $\Cal L_2$ in (2.5), the values at $(0,1)$ and $(0,-1)$ are the same
as the values for $\Cal L_1$, so (2.11) holds with the same values,
and Assumption 3.2 is satisfied. But not Assumption 3.1 since $\Cal
L_2$ is not strongly elliptic.

When $b_n\ne 0$, hence $\delta \ne 0$, neither of these symbols
satisfy the full  $\mu $-transmission condition
Definition 2.1 $2^\circ$, since second derivatives remove the $(i\pmb
b\cdot \xi )$-term so that the resulting symbol is even (with $\mu =a+\delta  $
replaced by $\mu =a$).

\endexample

Our answer to (1) is now the following (achieved in Section 5.4):

\proclaim{Theorem 2.4} Let $P$ satisfy Assumption {\rm 3.2}. For $\operatorname{Re}\mu -\frac12<t<
\operatorname{Re}\mu +\frac32$, $r^+P$  defines a continuous linear mapping
$$
r^+ P\colon  H^{\mu (t) }(\crnp)\to \ol H^{t -2a }(\rnp). \tag2.13
$$
\endproclaim

It is important to note that $r^+P$  then also makes good sense on
subsets of $H^{\mu (t) }(\crnp)$. In particular,  since $\E_\mu (\crnp)\cap
\E'(\rn)$ is a subset of  $H^{\mu (t)}(\crnp)$ for all $t$, the operator $r^+P$ is
well-defined on $\E_\mu (\crnp)\cap
\E'(\rn)$, mapping it into $ \bigcap_{t<\operatorname{Re}\mu +\frac32}\ol
H^{t-2a  }(\rnp)\subset \ol
H^{\operatorname{Re}\delta  +\frac32-a-\varepsilon  }(\rnp)$, any
$\varepsilon >0$, by (2.13). When  $\operatorname{Re}\delta
>-\frac12$ (always true under Assumption 3.1),  this is assured to be contained in $\ol
H^{1-a }(\rnp)$.

Our answer to (2) is (cf.\ Section 5.4):

\proclaim{Theorem 2.5} Let $P$ satisfy Assumption {\rm 3.1}. Then
$P=\widehat P+P'$, where  $P'$ is of order $2a-1$, and  $r^+\widehat P$ is a bijection from $H^{\mu
(t) }(\crnp)$ to $\ol H^{t -2a }(\rnp)$ for $\operatorname{Re}\mu -\frac12<t<\operatorname{Re}\mu
+\frac32$. In other words, there is unique solvability of {\rm
(2.7)} with $P$ replaced by $\widehat P$, in the mentioned spaces.

For $r^+P$ itself, there holds the regularity property:
Let $\operatorname{Re}\mu -\frac12<t<\operatorname{Re}\mu
+\frac32$, let $f\in
\ol H^{t-2a}(\rnp)$, and let  $u\in \dot H^{\sigma  }(\crnp) $ (for
some $\sigma >\operatorname{Re}\mu -\frac12$) solve the homogeneous
Dirichlet problem {\rm (2.7).}
Then $u\in H^{\mu (t)}(\crnp)$.
\endproclaim

The last statement shows a lifting of the regularity of $u$ in the elliptic
case, namely if it solves (2.7) lying in a low-order space $\dot H^{\sigma
}(\crnp) $, then it is in the  best possible $\mu $-transmission space
according to Theorem 2.4, mapping into the given range space $\ol H^{t-2a}(\crnp)$.
In other words, the {\it domain} of the homogeneous Dirichlet problem with range in
$\ol H^{t-2a}(\crnp)$ equals  $H^{\mu (t)}(\crnp)$.
 
The strategy for both theorems is, briefly expressed, as follows:
The first step is to replace $P=\Op(p(\xi ))$ by $\widehat P=\Op(\widehat p(\xi
))$, where $\widehat p(\xi )$ is better controlled at $\xi '=0$ and
$p'(\xi )=p(\xi )-\widehat p(\xi )$ is  $O(|\xi |^{2a-1})$ for $|\xi |\to \infty $.
The second step is to reduce $\widehat P$ to order 0 by composition
with "plus/minus order-reducing operators" $\Xi ^t_\pm=\Op((\ang{\xi
'}\pm i\xi _n)^t)$ ((3.11), (5.2)) geared to the value $\mu $ (recall $\mu '=2a-\mu $):
$$
\widehat Q=\Xi _-^{-\mu '}\widehat P \Xi _+^{-\mu } .\tag2.14
$$
Then the homogeneous symbol $q$ associated with $\widehat Q$ satisfies
the principal $0$-transmission condition.
The third step is to decompose $\widehat Q$ into a sum (when Assumption
3.2 holds) or a product (when Assumption 3.1 holds) of
operators whose action relative to the usual Sobolev spaces $\dot
H^s(\crnp)$ and $\ol H^s(\rnp)$ can be well understood, so that we can
show forward mapping properties and (in the strongly elliptic case) bijectiveness properties for
$\widehat Q$. The fourth step is to carry this over to forward mapping properties and (in the strongly elliptic case) bijectiveness properties for $\widehat P$.
The fifth and last step is to take $P'=P-\widehat P$ back into the picture and deduce the forward mapping resp.\ regularity
properties for the original operator
$P$.

It is the right-hand factor $\Xi _+^{-\mu }$ in (2.14) that is
the reason why the $\mu $-transmission spaces, defined by 
$$
H^{\mu
(t)}(\crnp)=\Xi _+^{-\mu }e^+\ol H^{t-\operatorname{Re}\mu }(\rnp),\tag2.15
$$
 enter.  Here $ e^+\ol
H^{t-\operatorname{Re}\mu  }(\rnp)$ has a jump at $x_n= 0$ when
$t>\operatorname{Re}\mu +\frac12$, and then the
 coefficient $x_n^\mu $ appears.
 
The analysis of $\widehat Q$ is based on  a Wiener-Hopf technique
(cf.\ Section 4) explained in
Eskin's book \cite{E81}, instead of the involvement of the extensive Boutet
de Monvel calculus used in \cite{G15}. 

An interesting feature of the results is that the $\mu $-transmission spaces have a
universal role, depending only on $\mu $ and not on the exact form of $P$.

Finally, we answer (3) by showing an integration by parts formula, based just on
Assumption 3.2.

\proclaim{Theorem 2.6} Let $P$ satisfy Assumption {\rm 3.2},
and assume moreover that $\operatorname{Re}\mu>-1 $,
$\operatorname{Re}\mu'>-1 $. For $u\in \E_\mu (\crnp)\cap \E'(\rn)$, $u'\in \E_{\bar\mu '} (\crnp)\cap \E'(\rn)$,
 there holds 
$$
\aligned 
\int_{\rnp}  Pu\,\partial_n\bar
u'\,dx&+\int_{\rnp}\partial_nu\,\overline{  P^*u'}\,dx\\
&=\Gamma (\mu +1){\Gamma(\mu '+1)}\int_{{\Bbb
R}^{n-1}}s_0\gamma _0(u/x_n^{\mu })\,{\gamma _0(\bar u'/x_n^{\mu '})}\,
dx',
\endaligned
\tag2.16
$$
where $s_0=e^{-i\pi \delta }p(0,1)$. The formula extends to $u\in H^{\mu (t)}(\crnp)$,  $u'\in H^{\bar\mu
'(t')}(\crnp)$, for  $t>\operatorname{Re}\mu +\frac12$,
$t'>\operatorname{Re}\mu '+\frac12$.
\endproclaim

The integrals over $\rnp$ in (2.16) are interpreted as dualities when
needed. The basic step in the proof is the treatment of one
order-reducing operator in Proposition 6.1, by an argument shown in detail in
\cite{G16, Th.\ 3.1, Rem. 3.2} and recalled in  \cite{G21, Th.\ 4.1}. 

In the proof of (2.16) in Section 6, the formula is first shown for the nicer operator
$\widehat P$, and thereafter extended to $P$.
(The formula (2.16) for $(-\Delta )^a$ in Ros-Oton and Serra
\cite{RS14, Th.\ 1.9} should have a minus sign on the boundary
contribution; this has been corrected by Ros-Oton in the survey \cite{R18, p.\ 350}.)

The theory will be carried further, to include ``large'' solutions of a
nonhomogeneous local Dirichlet problem, and to show
regularity results and a ``halfways Green's formula'', see Section 6,
but we shall leave those aspects out of this preview.

The example $L_1$ in (2.3) is a special case of the operator
$L=\Op(\Cal L(\xi ))$, where
$\Cal L(\xi )=\Cal A(\xi
)+i\Cal B(\xi )$ with $\Cal A(\xi )$ real, {\it  even} in $\xi $ and positive, and $\Cal B(\xi
)$ real and {\it odd} in $\xi $. There are more details below in  (3.5)ff.\ (this
stands for (3.5) and the near following text) and Examples
5.9, 6.5, 7.4. $L$ was first studied in \cite{DRSV21} (under certain restrictions on $\mu $), and our results apply to
it. Theorem 2.6 gives an alternative proof for the same integration by
parts formula,
established in \cite{DRSV21, Prop.\ 1.4} by extensive real function-theoretic
methods.

The result on the integral over $\rnp$ is combined in \cite{DRSV21} with localization
techniques to get an interesting result for curved domains, and it is
our hope that the present results for more general strongly elliptic
operators can be used in a similar way.

\subhead 3. The principal $\mu $-transmission condition \endsubhead

\subsubhead 3.1 Analysis of homogeneous symbols \endsubsubhead

Let $p(\xi )$ be a complex function on $\rn$ that is homogeneous of degree $m$ in $\xi $,
and let $\nu \in \rn$ be a unit vector. For a complex number $\mu $, we shall say
that $p$ satisfies {\it the principal $\mu $-transmission condition in the
direction $\nu $}, when 
$$ 
p(-\nu )=e^{i\pi (m-2\mu )}p(\nu ).\tag3.1
$$
When $p(\nu )\ne 0$, we can rewrite (3.1) as
$$
e^{i\pi (m-2\mu )}=\tfrac {p(-\nu )}{p(\nu )}, \text{ i.e., }\mu
=\tfrac m2-\tfrac 1{2\pi i}\log \tfrac {p(-\nu )}{p(\nu )},
$$
where $\log$ is a complex logaritm. This determines the possible $\mu
$ up to addition of an integer.

The (full) $\mu $-transmission property defined in \cite{G15} demands much
more, namely that
$$
\partial_\xi ^\alpha p(-\nu )=e^{i\pi (m-2\mu -|\alpha |)}\partial_\xi
^\alpha p(\nu ),\text{ all }\alpha \in \N_0^n.\tag3.2
$$
Besides assuming infinite differentiability, this is a stronger
condition than (3.1) in particular 
because of the
requirements it puts on derivatives of $p$ transversal to $\nu $.

To analyse this we observe that when a (sufficiently smooth) function $f(t)$ on $\R\setminus \{0\}$ is
homogeneous of degree $m\in\R$, then it 
has the form, for some $c_1,c_2\in\C$,
$$
f(t)=\cases c_1t^m\text{ for }t>0,\\
c_2(-t)^m\text{ for }t<0,\endcases
$$
  and its derivative outside $t=0$ is a function homogeneous of degree $m-1$
  satisfying
$$
\partial_tf(t)=\cases c_1mt^{m-1}\text{ for }t>0,\\
-c_2m(-t)^{m-1}\text{ for }t<0.\endcases
$$ In particular, if $c_1\ne 0$, $m\ne 0$,
$$
f(-1)/f(1)=c_2/c_1,\quad \partial_tf(-1)/\partial_tf(1)=-c_2/c_1.
$$
In the case $m=0$, $f$ is constant for $t>0$ and $t<0$, and the
derivative is zero there.

Thus, when $p(\xi )$ is a (sufficiently smooth) function on $\rn\setminus\{0\}$ that is
homogeneous of degree $m \ne 0$, and we consider it on a two-sided ray $\{t\nu
\mid t\in\R\}$ where $\nu $  is a unit vector and $p(\nu )\ne 0$, then
$$
p(-\nu )=c_0p(\nu )\implies \partial_tp(t\nu )|_{t=-1}=-c_0\partial_tp(t\nu 
)|_{t=1}.\tag 3.3
$$

So for example, when $\nu $ is the inward normal $(0,1)=\{(\xi ',\xi
_n)\mid \xi '=0,\xi _n=1\}$ to $\rnp$,  
$$
p(0,-1 )=c_0p(0,1 )\implies \partial_{\xi _n}p(0,-1 )=-c_0\partial_{\xi _n}p(0,1).
$$

For $p(\xi )$ satisfying (3.1), this means that 
when 
$p(\nu )\ne 0$, it will also satisfy 
$$
\partial_tp(t\nu )|_{t=-1}=e^{i\pi (m-2\mu -1)}\partial_tp(t\nu )|_{t=1},
$$ 
in view of (3.3). This argument can be repeated, showing that 
$$
\partial_t^kp(t\nu )|_{t=-1}=e^{i\pi (m-2\mu -k)}\partial_t^kp(t\nu 
)|_{t=1},\tag3.4
$$
for all $k\in \N$ (possibly vanishing from a certain step on). On the other hand, we cannot infer that other
derivatives $\partial_\xi ^\alpha $ of $p$ have the
property (3.2); an example will be given below.

In general, $\mu $ takes different values for
different $\nu $. When $\Omega $ is a sufficiently smooth subset of
$\rn$ with interior normal $\nu (x)$ at boundary points
$x\in\partial\Omega $, we say that $p$ satisfies {\it the principal $\mu
$-transmission condition at} $\Omega $ if $\mu (x)$ is a function on
$\partial\Omega $ such that (3.1) holds with this $\mu (x)$ at boundary
points $x\in\partial\Omega $. For $\Omega =\rnp$, the normal $\nu $ equals
$(0,1)$ at all boundary points and  $\mu $ is a constant; this is the
situation considered in the present paper.

In \cite{G21} we have studied a special class of symbols first
considered by Dipierro, Ros-Oton, Serra and Valdinoci in
\cite{DRSV21}:
$$
\Cal L(\xi )=\Cal A(\xi
)+i\Cal B(\xi ),\tag3.5 
$$
the functions being $C^\infty $ for $\xi \ne 0$ and homogeneous in $\xi
$ of degree $2a>0$ ($a<1$), and
 where
 $\Cal A(\xi )$ is real and  even in $\xi $ (i.e., $\Cal A(-\xi )=\Cal A(\xi )$), $\Cal B(\xi
)$ is real and odd in $\xi $ (i.e., $\Cal B(-\xi )=-\Cal B(\xi )$),
 and $\Cal L$ is strongly elliptic (i.e.,  $\Cal A(\xi )>0$ for $\xi \ne
0$). As shown in \cite{G21, Sect.\ 2}, $\Cal L$ satisfies (3.1) on each unit
vector  $\nu $, for $m=2a$ and
$$
\mu(\nu ) =a+\delta (\nu ),\text{ with }\delta (\nu )=\tfrac1\pi 
\operatorname{Arctan}b, \; b=\Cal B(\nu )/\Cal A(\nu );\tag3.6
$$
this follows straightforwardly (as in Example 2.3) from the observation that $\Cal L(-\nu )/\Cal L(\nu
)=(1-ib)/(1+ib)$, $b= B(\nu )/\Cal A(\nu )$. It then also satisfies
(3.4) with this $\mu $.

But the full $\mu $-transmission condition need not hold.
For example, the
symbol $\Cal L_1(\xi )=|\xi |+i\pmb b\cdot \xi $ in (2.3) (with $\pmb b\in \rn$)
satisfies the principal $\mu $-transmission condition for $\nu =(0,1)$
with $\mu
=\frac12+\delta $, $\delta\ne 0$ if $b_n\ne 0$, whereas
$$
\partial_{\xi _1}^2\Cal L_1=|\xi '|^2/|\xi |^3
$$
and its derivatives satisfy the conditions in (3.2) for $\nu =(0,1)$ with $\mu $ replaced by $a$.

The statement in \cite{G21, Th.\ 3.1} that solutions of the
homogeneous Dirichlet problem have a structure with the factor $x_n^\mu
$, was quoted from \cite{G15} based on the full $\mu $-transmision condition, and therefore
applies to $L=\Op(\Cal L)$ 
when $\Cal B=0$ (a case belonging to \cite{G15}), but not in general when $\Cal B \ne 0$. Likewise, the
integration by parts formulas for $L$ derived in \cite{G21} using details from the Boutet de Monvel
calculus are justified when $\Cal B=0$ or when other operators $P$
satisfying the full $\mu $-transmission condition are inserted, but
not in general when $\Cal B\ne 0$. Fortunately, there are  cruder
methods that do lead to such results, on the basis of the principal
$\mu $-transmission condition alone, and that is what we show in this
paper.

The treatment of $\Cal L$ will be incorporated in a treatment of general strongly elliptic
homogeneous symbols in the following. This requires that we allow complex values of $\mu $.

Let $P=\Op(p(\xi ))$ be defined by (2.1)
from a symbol $p(\xi )$ that is $C^1 $ for $\xi \ne 0$, homogeneous of
order $m=2a>0$, and now also strongly elliptic (2.4).
To fix the ideas, we shall consider the operator relative to
the set $\rnp$, with interior normal $\nu =(0,1)$.
Denote $p(\xi )|\xi |^{-2a}=p_1(\xi )$; it is homogeneous
of degree 0. Both $p$ and $p_1$ take values in a closed subsector of
$\{z\in\C\mid \operatorname{Re}\xi
_n>0\}\cup \{0\}$. For any $\xi '\in \R^{n-1}$, one has for $+1$ and
$-1$ respectively,
$$
\lim_{\xi _n\to \pm\infty }p_1(\xi ',\xi _n)=\lim_{\xi _n\to \pm\infty
}p_1(\xi '/|\xi _n|,\pm 1)=p_1(0,\pm 1)=p(0,\pm 1).
$$
With the logarithm $\log z$ defined to be positive for real $z>1$,
with a cut along the negative real axis, denote $\log p(0,\pm 1)=\alpha _\pm$; here
$\operatorname{Re}\alpha _\pm=\log|p(0,\pm 1)|$ and
$\operatorname{Im}\alpha _\pm$ is the argument of $p(0,\pm 1)$. With this notation,
$$
p(0,-1)/p(0,1)=e^{\alpha _-}/e^{\alpha _+}=e^{\alpha _--\alpha _+},
$$
 so (3.1) for $m=2a$ holds with $\nu =(0,1)$ when $\alpha _--\alpha _+=i\pi (2a-2\mu )$, i.e., 
$$
\mu =a+\delta \text{ with }\delta =(\alpha _+-\alpha _-)/2\pi i;\tag3.7
$$
this $\mu $ is the {\it factorization index}.
These calculations were given in \cite{G15, Sect.\ 3} (with $m=2a$), and are in
principle consistent
with the determination of the factorization index by Eskin in  \cite{E81, Ex.\ 6.1} (which
has different plus/minus conventions because of a different definition
of the Fourier transform).

Since $p(\xi )$ takes values in $\{\operatorname{Re}z>0\}$ for $\xi
\ne 0$, both $p(0,1)$ and $p(0,-1)$ lie there and the difference between their arguments is
less than $\pi $, so $|\operatorname{Im}(\alpha _+-\alpha _-)/2\pi
|<\frac12$; in other words $$
|\operatorname{Re}\delta |<\tfrac12.\tag3.8
$$
Note that $\delta $ is real in the case (3.5).

We collect the information on $P$ in the following description:

\proclaim{Assumption 3.1} The operator $P=\Op(p(\xi ))$ is defined
from a symbol $p(\xi )$ that is $C^1 $ for $\xi \ne 0$, homogeneous of
order $m=2a>0$, and strongly elliptic {\rm (2.4)}.
It satisfies the principal $\mu $-transmission condition in the
direction $(0,1)$:
$$
p(0,-1)=e^{i\pi (m-2\mu )}p(0,1),
$$
with
$\mu $ equal to the factorization index  $\mu =a+\delta $ derived
around {\rm
(3.7)}, and  $|\operatorname{Re}\delta |<\tfrac12$. Denote $\mu '=2a-\mu =a-\delta $.
\endproclaim

In the book \cite{E81}, the case of constant-coefficient
pseudodifferential operators considered on $\rnp$ is
studied in $\S\S$4--17, and the calculations rely on the principal transmission condition up to and including $\S$9. From $\S$10 on, additional
conditions on transversal
derivatives are required (the  symbol class $D^{(0)}_{\alpha
+i\beta }$ seems to correspond to our full 0-transmission condition,
giving operators preserving smoothness up to the boundary). In the
following, we draw on some of the points  made in $\S\S$6--7 there.

For an operator $A$ defined from a homogeneous symbol $a(\xi )$, the
behavior at zero can be problematic to deal with. In \cite{E81,
\S7} there is introduced a
technique that leads to a nicer operator, in the context of operators
relative to $\rnp$: One eliminates the
singularity at $\xi '=0$ by
replacing the homogeneous symbol $a(\xi ',\xi _n)$ by $$
\widehat a(\xi ',\xi
_n)=a(\ang{\xi '}\xi '/|\xi '|,\xi _n),\tag3.9
$$
the corresponding operator denoted $\widehat A$.
 (In comparison with \cite{E81} we have
 replaced the factor $1+|\xi '|$ used there by $\ang{\xi '}=(1+|\xi '|^2)^{\frac12}$.)
It is
shown there that when $a(\xi )$ is homogeneous of degree $\alpha +i\beta $, then
$$
a'(\xi )=a(\xi )-\widehat a(\xi ) \text{ is }O(|\xi |^{\alpha -1})\text{
for }|\xi |\ge 2
,\tag3.10
$$
hence is of lower order in a certain sense.
Many results with Sobolev estimates are then shown primarily for the
``hatted'' version $\widehat A=\Op(\widehat a)$, and supplied afterwards with information on $A'=\Op(a')$. Indeed, we shall see that the results we are
after for our operators $P=\Op(p)$, can be obtained in a manageable way for
$\widehat P= \Op(\widehat p)$, and then extended to $P$ by a supplementing analysis of $P'$. The important thing is that special properties
with respect to $\xi _n$, such as holomorphic extendability into
$\C_+$ or $\C_-$, are not disturbed when $a$ is replaced by $\widehat a$.

Some of the results that we shall show do not require ellipticity of
$P$. We therefore introduce also a weaker assumption:

\proclaim{Assumption 3.2} The operator $P=\Op(p(\xi ))$ is defined
from a symbol $p(\xi )$ that is $C^1 $ for $\xi \ne 0$, homogeneous of
order $m=2a>0$, and satisfies the principal $\mu $-transmission
condition in the direction $(0,1)$
with $\mu =a+\delta $
for some $\delta \in\C$. Denote $a-\delta =\mu '$.
\endproclaim

For
the symbols $p$ considered in the rest of the paper, we assume at least that
Assumption 3.2 holds. As noted earlier, when $P$ satisfies (3.1) for some $\mu
$, it also does so with $\mu $ replaced by $\mu +k$, $k\in\Z$. The
precision in Assumption 3.1, that $\mu $ should equal the factorization
index, is needed for elliptic solvability statements.

\subsubhead 3.2 Reduction to symbols of order $0$ \endsubsubhead

Consider the symbols of ``order-reducing'' operators (more on them in
Section 4):
$$
\aligned
\chi _{0,\pm}^t(\xi )&=(|\xi '|\pm i \xi _n)^t;\text{ consequently }\\
\widehat \chi _{0,\pm}^t(\xi
)&=
(|\ang{\xi '}\xi '/|\xi '||\pm i \xi _n)^t=(\ang{\xi '}\pm i\xi _n)^t=\chi
_\pm^t(\xi );
\endaligned\tag 3.11
$$
the last entry is the usual notation.
Together with our  symbol $p(\xi )$ of order $2a$, we shall
consider its reduction to a symbol $q$ of order 0 defined by:
$$
q(\xi )=\chi  _{0,-}^{-\mu '}p(\xi )\chi  _{0,+}^{-\mu },\text{ hereby
}p(\xi )=\chi  _{0,-}^{\mu '}q(\xi )\chi  _{0,+}^{\mu }.\tag3.12
$$
 The ``hatted'' version is:
$$
\widehat q(\xi )= \chi  _-^{-\mu '}\widehat{p}(\xi )\chi
_+^{-\mu },\text{ hereby } \widehat{p}(\xi )=\chi  _-^{\mu '}
\widehat q(\xi )\chi  _+^{\mu }.\tag3.13
$$
Here $q$ is continuous and  homogeneous of degree 0 for $\xi \ne 0$; it
is  $C^1 $ in $\xi _n$ there, and $C^1 $ in $\xi '$ for $\xi '\ne 0$ with bounded
first derivatives on $|\xi |=1$. Since $i=e^{i\pi /2}$,
$$
\aligned
q(0,1)&=(-i)^{\mu -2a }p(0,1)i^{-\mu  }=i^{2a-2\mu 
}p(0,1)=e^{i\pi (a-\mu ) }p(0,1),\\
q(0,-1)&=(+i)^{\mu -2a }p(0,-1)(-i)^{-\mu  }=i^{2\mu -2a
}e^{i\pi (2a-2\mu )}p(0,1)\\
&=e^{i\pi (a-\mu ) }p(0,1)=q(0,1),
\endaligned
$$
so $q$ satisfies the principal 0-transmission condition in the
direction $\nu =(0,1)$:
$$
q(0,-1)=q(0,1).\tag3.14
$$
In view of (3.1)--(3.4), we have moreover when $p(0,1)\ne 0$ that
$$
 \partial_{\xi _n}q(0,-1)=-\partial_{\xi _n}q(0,1).\tag3.15
$$
Note that since $\mu -a=\delta $, $q(0,1)=e^{-i\pi \delta }p(0,1)$.
We shall
denote$$
s_0=q(0,1)=e^{-i\pi \delta }p(0,1).\tag3.16
$$

In the case $p=\Cal L$ in (3.5)--(3.6), $\Cal L(0,1)=e^{i\pi \delta }|\Cal
L(0,1)|$ with $\delta $ real, so 
$$
s_0=e^{-i\pi \delta }\Cal L(0,1)=|\Cal L(0,1)|=(\Cal A(0,1)^2+\Cal B(0,1)^2)^{\frac12}\text{ then.}\tag3.17
$$

\subhead 4. The Wiener-Hopf decomposition \endsubhead

\subsubhead 4.1 The sum decomposition \endsubsubhead

Since $p(\xi )$ is only assumed to satisfy the principal $\mu
$-transmission condition, $q(\xi )$ will in general only satisfy the
principal $0$-transmission condition, not the full one, so the
techniques of the Boutet de Monvel calculus brought forward in \cite{G15}
are not available.
Instead we go back to a more elementary application of the original
Wiener-Hopf method \cite{WH31}.

When 
$b(\xi _n)$ is a function on $\R$, denote
$$
\aligned
b_+(\xi _n+i\tau )=\frac{i}{2\pi }\int_{\R}\frac {b(\eta _n)}{\eta
_n-\xi _n-i\tau }\,d\eta _n\text{ for }\tau < 0,
\\
b_-(\xi _n+i\tau )=\frac{-i}{2\pi }\int_{\R}\frac {b(\eta _n)}{\eta
_n-\xi _n-i\tau } \,d\eta _n\text{ for }\tau > 0,
\endaligned \tag4.1
$$
when the integrals have a sense. When $b$ is suitably nice, $b_+$ is
holomorphic in $\xi _n+i\tau $ for $\tau <0$ and extends to a
continuous function on $\overline{\C}_-$ (also denoted $b_+$), $b _-$ has these
properties relative to $\overline{\C}_+$, and $b(\xi _n)=b_+(\xi
_n)+b_-(\xi _n)$ on $\R$. With the notation of spaces $H$, $H^\pm$
introduced by Boutet de Monvel in \cite{B71}, denoted $\Cal H$, $\Cal H^\pm$ in our
subsequent works,
the decomposition holds for $b\in \Cal H$ with $b_{\pm}\in \Cal H^\pm$
on $\R$. Since we are presently dealing with functions with  cruder
properties, we shall instead apply a useful lemma shown in \cite{E81, Lemma
6.1}:

\proclaim{Lemma 4.1} Suppose that $b(\xi ',\xi _n)$ is homogeneous of
degree $0$ in $\xi $, is $C^1$ for $\xi '\ne 0$, and satisfies
$$
|b(\xi ',\xi _n)|\le C|\xi '|\,|\xi |^{-1},\quad |\partial_jb(\xi
 ',\xi _n)|\le C|\xi |^{-1}\text{ for }j\le n-1.\tag4.2 
$$
Then the function defined for $\tau <0$ by
$$
b_+(\xi ',\xi _n+i\tau )=\frac{i}{2\pi }\int_{\R}\frac {b(\xi ',\eta _n)}{\eta
_n-\xi _n-i\tau }\,d\eta _n\tag4.3
$$
 is holomorphic with respect to $\xi _n+i\tau $ in $\C_-$, is homogeneous
 of degree $0$, extends by continuity with respect to $(\xi ',\xi _n+  i\tau )\in \overline{\C}_-$ for $|\xi
 |+|\tau |>0$, $\tau \le 0$, and satisfies the estimate
$$
|b_+(\xi ',\xi _n+i\tau )|\le C_\varepsilon  |\xi '|^{1-\varepsilon 
 }(|\xi |+|\tau |)^{\varepsilon  -1}, \text{ any }\varepsilon
 >0.\tag4.4 
$$

 There is an analogous statement for $b_-$ with $\C_-$ replaced by $\C_+$.
\endproclaim

The symbol $q$ derived from $p$ by (3.12)ff.\ satisfies
$$
q(\xi )=s_0+f(\xi ),
$$
where $f$ is likewise homogeneous of degree 0, and has $f(0,1)=f(0,-1)=0$.
We make two applications of Lemma 4.1. One is, under Assumption 3.2,  to apply it directly to $f$ to get a sum decomposition
$f=f_++f_-$ where the terms extend holomorphically to $\C_-$ resp.\
$\C_+$ with respect to $\xi _n$; this will be convenient in
establishing the forward mapping properties and integration by parts formula  
for the present operators. The other is, under Assumption 3.1, to apply the lemma to the function $b(\xi )=\log
q(\xi )$ to get a sum decomposition of $b$ and hence a {\it
factorization} of $q$; this is used to show that $P$ has appropriate
solvability  properties (the solutions exhibiting a singularity $x_n^\mu $ at the
boundary). 

We show that $f$ has the properties required for Lemma 4.1 as follows:
To see that (4.2) is verified by $f$, note that the second inequality
follows since $\partial_jf$ is bounded on the unit sphere $\{|\xi
|=1\}$ and homogeneous of degree $-1$. For the first inequality we
have, when $\xi _n>|\xi '|$ (hence $|\xi '/\xi _n|< 1$),
$$
|f(\xi ',\xi _n)|=\Bigl|q\Bigl(\frac{\xi '}{\xi _n},1\Bigr)-q(0,1)\Bigr|
\le
 {\sum}_{j<n}\Bigl|\frac{\xi _j}{\xi _n}\Bigr|\sup_{|\eta '|\le
 1}|\partial_jq(\eta ',1)|\le C \frac{|\xi '|}{|\xi _n|}\le C'
 \frac{|\xi '|}{|\xi |},
\tag4.5 
$$
 using the mean value theorem and the fact that $|\xi _n|\sim |\xi |$
 when $|\xi _n|\ge |\xi '|$. A similar estimate is found for $\xi _n<
 -|\xi '|$.
 For $|\xi _n|\le |\xi '|$, we use that $q$ is bounded, so that
 $|q(\xi )-s_0||\xi |/|\xi '|\le c |q(\xi )-s_0||\xi '|/|\xi '|\le c'$.
We have obtained:

\proclaim{Proposition 4.2} When $p$ satisfies Assumption {\rm 3.2} and $q$ is derived from $p$ by {\rm (3.12)}ff.\,
then there is a
sum decomposition of $f=q-s_0$:
$$
q(\xi )-s_0=f_+(\xi )+f_-(\xi ), 
$$
where $f_+(\xi ',\xi _n)$  is holomorphic with respect to $\xi _n+i\tau $ in $\C_-$, and
 continuous with respect to $(\xi ',\xi _n+  i\tau )\in \overline{\C}_-$ for $|\xi
 |+|\tau |>0$, $\tau \le 0$, and satisfies  estimates
$$
|f_+(\xi ',\xi _n+i\tau )|\le C_\varepsilon  |\xi '|^{1-\varepsilon 
 }(|\xi |+|\tau |)^{\varepsilon  -1}, \text{ any }\varepsilon
 >0,\tag4.6 
$$
and $f_-$
has the analogous properties with $\C_-$ replaced by $\C_+$.
\endproclaim

For the corresponding hatted symbol, we then have 
 $\widehat q=s_0+\widehat f_+ +\widehat f_-$, with $\widehat f_\pm$ defined from $ f_\pm$.
 They have similar holomorphy properties, and satisfy
 estimates as in (4.6) with $|\xi '|$ replaced by $\ang{\xi '}$.

\subsubhead 4.2 The product decomposition \endsubsubhead

In order to obtain a factorization for symbols satisfying Assumption
3.1, we shall study $\log q$. By
the strong ellipticity, $q(\xi )\ne 0$
for $\xi \ne 0$. Moreover,  $p (\xi )|\xi
|^{-2a}=\chi _{0,-}^{-a}p(\xi )\chi _{0,+}^{-a}$ takes values in a
subsector of $\{z\in\C\mid\operatorname{Re}z>0\}$ and the multiplication by
$\chi _{0,-}^{\delta }$ and $\chi _{0,+}^{-\delta }$ gives the
function $q$
taking  values in the sector $\{z\in\C\mid |\arg z|\le \pi(\frac12+|\operatorname{Re}\delta |)\}$ disjoint
from the negative real axis. So the logarithm is well-defined with
inverse exp.

Assume first that $s_0=1$; this can simply be
obtained by dividing out $q(0,1)$. The function $b(\xi )=\log q(\xi )$
is homogeneous of degree 0 and has $b(0,1)=b(0,-1)=0$ and
 the appropriate continuity properties, and bounds on first
derivatives, so the same proof as for $f$ applies to $b$ to give the
decomposition $b=b_++b_-$.  Then we define $q^\pm=\exp (b_\pm)$, they are
homogeneous of degree 0. For
example, 
$$
q^+=1+g^+, \text{ where }g^+={\sum}_{k\ge 1}(b_+)^k.
$$
Here $|b_+(\xi )|\le C_\varepsilon |\xi '|^{1-\varepsilon }|\xi |^{-1+\varepsilon }$, and 
there is a constant $C'_\varepsilon $ such that   $C_\varepsilon |\xi
'|^{1-\varepsilon }|\xi |^{-1+\varepsilon }\le \frac12$ for $|\xi _n|\ge
C'_\varepsilon |\xi '|$. On this set the series for $g^+$ converges with $|g^+|\le
|b_+|$, hence $g^+$ satisfies an estimate of the form (4.4) there. It
likewise does so on the set  $|\xi _n|\le
C'_\varepsilon |\xi '|$ since $|\xi '|\sim |\xi |$ there. There are
similar results for $q^-=\exp (b_-)=1+g^-$ with $\C_-$ replaced by $\C_+$.
This shows:

\proclaim{Proposition 4.3}   When $p$ satisfies Assumption {\rm 3.1} and $q$ is derived from $p$ by {\rm (3.12)}ff.\
and satisfies $s_0=1$, then there is a
factorization of $q$:
$$
q(\xi )=q^-(\xi )q^+(\xi ), 
$$
where $q^+(\xi ',\xi _n)$  is holomorphic with respect to $\xi _n+i\tau $ in $\C_-$, and
 continuous with respect to $(\xi ',\xi _n+  i\tau )\in \overline{\C}_-$ for $|\xi
 |+|\tau |>0$, $\tau \le 0$. Moreover, $g^+=q^+-1$  satisfies  estimates
$$
|g^+(\xi ',\xi _n+i\tau )|\le C_\varepsilon  |\xi '|^{1-\varepsilon 
 }(|\xi |+|\tau |)^{\varepsilon  -1}, \text{ all }\varepsilon
 >0,\tag4.7 
$$
and $q^-$, $g^-=q^--1$
have the analogous properties with $\C_-$ replaced by $\C_+$. The symbols
are homogeneous of degree $0$, and $q^+$ and $q^-$
are elliptic.

For general $s_0$, we apply the factorization to $q_0=s_0^{-1}q$, so
that $q_0=q_0^-q_0^+$; then $q=q^-q^+$ with $q^-=s_0q_0^-=s_0(1+g^-)$ and $q^+=q_0^+=1+g^+$.
\endproclaim

The ellipticity follows from the construction as $\exp({b_\pm})$, or one
can observe that the product $q^+q^-=q$ is elliptic (i.e.,
nonzero for $\xi \ne 0$).

The notation with upper index $\pm$ is chosen here to avoid confusion with
the lower $+$ used later to indicate truncation, $P_+=r^+Pe^+$.

Turning to the corresponding hatted symbols, we have obtained
 $\widehat q=\widehat q^- \widehat q^+$, with $\widehat q^\pm$,
 $\widehat g^\pm$ defined
 from $ q^\pm$, $g^\pm$, respectively.
  They have similar holomorphy properties,  the $\widehat q^\pm$
 are elliptic, and the $\widehat g^\pm$ satisfy
 estimates as in (4.7) with $|\xi '|$ replaced by $\ang{\xi '}$:
$$
|\widehat g^+(\xi ',\xi _n+i\tau )|\le C_\varepsilon  \ang{\xi '}^{1-\varepsilon 
 }(\ang{\xi }+|\tau |)^{\varepsilon  -1}, \text{ all }\varepsilon
 >0.\tag4.8 
$$

\subhead 5. Mapping properties and the homogeneous Dirichlet problem \endsubhead

\subsubhead 5.1 Some function spaces \endsubsubhead

First recall some terminology: $\Cal E'(\rn)$ is the space of
distributions on $\rn$ with
compact support, $\Cal S(\rn)$ is the Schwartz space of $C^\infty
$-functions $f$ on $\rn$ such that $x^\beta D ^\alpha f$ is bounded
for all $\alpha ,\beta $, and $\Cal S'(\rn)$ is its dual space of
temperate distributions. $\ang\xi $ stands for $(1+|\xi
|^2)^\frac12$. We denote by $r^+$ the operator restricting
distributions on $\rn$ to distributions on $\rnp$, and by $e^+$ the
operator extending functions on $\rnp$ by zero on
$\rn\setminus\rnp$. Then $r^+\Cal S(\rn)$ is denoted $\Cal
S(\crnp)$. 
The following notation for $L_2$-Sobolev spaces will be used, for $s\in\R$:
$$
\aligned
H^s (\Bbb R^n)
&=\{u\in \Cal S'(\Bbb R^n)\mid \langle{\xi }\rangle^s\F
u\in L_2 (\Bbb R^n)\},\\
\overline H ^s(\rnp )&=r^+H^s (\Bbb R^n), \text{ the {\it restricted} space},\\
\dot H^s (\crnp )&=\{u\in H^s (\Bbb R^n)\mid
\operatorname{supp}u\subset \crnp \}, \text{ the  {\it supported} space},
\endaligned\tag5.1
$$
as in our earlier papers on fractional-order operators. An elaborate
presentation of $L_p$-based spaces was given in \cite{G15}. (The notation
with dots and overlines stems from H\"ormander \cite{H85, App.\ B.2} and
is practical in formulas where both types of spaces occur. There are other
notations without the overline, and where the dot is replaced by a
ring or twiddle.)

Here $ \ol H^s(\rnp )$ identifies with the dual space of
$\dot H^{-s}(\crnp )$ for all $s\in\R$ (the duality extending the
$L_2(\rnp)$ scalar product). When $|s|<\frac12$, there is an identification of $\dot
H^s(\crnp )$ with $\ol H^s(\rnp )$ (more precisely with $e^+\ol H^s(\rnp )$).
 The trace operator $\gamma _0\colon u\mapsto \lim_{x_n\to 0+}u(x',x_n)$ extends
 to a continuous mapping $\gamma _0\colon \ol H^s(\rnp)\to
 H^{s-\frac12}(\R^{n-1})$ for $s>\frac12$.

The order-reducing operators $\Xi _\pm^t$ are defined for $t\in{\Bbb
C}$ by  $\Xi _\pm^t=\Op(\chi _\pm^t)$, where $\chi _\pm^t=(\ang{\xi '}\pm i\xi _n)^t$, cf.\ (3.11).
These operators have the
homeomorphism properties: 
$$
\Xi ^{t }_+\colon \dot H^s(\crnp )\simto
\dot H^{s- \operatorname{Re}t }(\crnp),\quad
r^+\Xi ^{t }_{-}e^+\colon \ol H^s(\rnp )\simto
\ol H^{s- \operatorname{Re}t } (\rnp ),\text{ all }s\in\R,
t\in\C;\tag5.2
$$
$r^+\Xi
 ^{ t }_{-}e^+$ is often denoted $\Xi ^{ t}_{-,+}$ for short. For each $t\in\C$, the operators $\Xi ^t _{+}$ and $\Xi ^{\overline t}_{-,+}$  identify with each other's adjoints
over $\crnp$ (more comments on this in \cite{G15, Rem.\ 1.1}). Recall also the simple composition rules (as noted e.g.\
 in \cite{GK93, Th.\ 1.2}):
$$
\Xi ^s_+\Xi ^t_+=\Xi
^{s+t}_+,\quad \Xi ^s_{-,+}\Xi ^t_{-,+}=\Xi
^{s+t}_{-,+}\text{ for $s,t\in\C$. }
$$

We define
$$
\E_\mu (\crnp)=e^+x_n^\mu C^\infty (\crnp)\text{ when
}\operatorname{Re}\mu >-1, \tag5.3
$$
 and $\E_\mu (\crnp)$
is defined successively as the linear hull of first-order derivatives
of elements of $\E_{\mu +1}(\crnp)$ when $\operatorname{Re}\mu \le -1$
(then distributions supported in the boundary can occur).
The spaces were introduced in H\"ormander's unpublished lecture notes
\cite{H66} and are presented in \cite{G15} (and with a different
notation in \cite{H85, Sect.\ 18.2}), and they satisfy for all $\mu $
(cf.\ \cite{G15, Props.\ 1.7, 4.1}):
$$
\E_\mu (\crnp)\cap \Cal E'(\rn)\subset\Xi _+^{-\mu }e^+ {\bigcap}_s\ol H^s(\crnp).\tag5.4
$$
A sharper statement follows from \cite{G21, Lemma 6.1} (when $\operatorname{Re}\mu >-1$):
$$
e^+x_n^{\mu }\Cal S(\crnp)= \Xi
_+^{-\mu }e^+\Cal S(\crnp).
\tag5.5
$$

\subsubhead 5.2 Mapping properties of the zero-order operator
$\widehat Q$ in Sobolev spaces \endsubsubhead

Let $P$ satisfy Assumption 3.1, and consider $\widehat
Q^{\pm}=\Op(\widehat q^{\pm})$, defined from the symbols
$q^{\pm}(\xi )$ introduced in Proposition 4.3. Since $\widehat q^\pm$ are bounded
symbols with bounded inverses, and extend holomorphically in $\xi _n$
into $\C_-$ resp.\ $\C_+$,  
$$
\widehat Q^+\colon \dot H^s(\crnp)\simto  \dot H^s(\crnp)\text{ and
 }\widehat Q^-_+=r^+\widehat Q^-e^+\colon \ol H^s(\rnp)\simto  \ol H^s(\rnp), \text{
 for all }s\in
 \R; \tag5.6
 $$ the latter follows since $r^+\widehat Q^-e^+$ is the adjoint of
 $\Op(\overline{\widehat q^-})$ over $\rnp$, where $\Op(\overline{\widehat q^-})$
 defines homeomorphisms in $\dot H^s(\crnp)$ (since $\overline{\widehat q^-}$
 has similar properties as $\widehat q^+$). The inverses $(\widehat Q^{\pm})^{-1}=\Op((\widehat q^\pm)^{\,-1})$ have
 similar homeomorphism properties. Since $\ol H^s (\rnp)=\dot H^s (\crnp)$
 for $|s|<\frac12$, it follows that we also have for $|s|<\tfrac12$:
 $$
 \widehat Q^+_+=r^+\widehat Q^+e^+\colon \ol H^s (\rnp)\simto \ol H^s
 (\rnp),\quad\widehat Q^-_+\widehat Q^+_+\colon \ol H^s (\rnp)\simto \ol H^s
 (\rnp). 
 $$

If $q$ satisfies the full 0-transmission condition, we are in the case
studied in \cite{G15}, and the bijectiveness in $\ol H^s(\rnp)$ can
be lifted to all higher $s$ by use of elements of the Boutet de Monvel
calculus, as accounted for in the proof of \cite{G15, Th.\ 4.4}. The
symbol $q$
presently considered is only known to satisfy the principal 
0-transmission condition (and possibly a few more identities).
We shall here show that a lifting is possible in general up to
$s<\frac32$.

\proclaim{Proposition 5.1} Let $P$ satisfy Assumption {\rm 3.1}, and
consider $\widehat Q^+=\Op(\widehat q^+)$ derived from it in Section {\rm 3}.

For any  $-\frac12<s<\frac32$, $\widehat Q^+_+=r^+\widehat Q^+e^+$ is continuous
$$
r^+\widehat Q^+e^+\colon \ol H^{s}(\rnp)\to \ol H^{s }(\rnp),\tag5.7
$$
and the same holds for the operator  $((\widehat Q^+)^{-1})_+$ defined
from its inverse $(\widehat Q^+)^{-1}$.

In fact, {\rm (5.7)} is a homeomorphism, and the inverse of $\widehat
Q^+_+$ is  $((\widehat Q^+)^{-1})_+$.
\endproclaim

\demo{Proof}  
We already have the mapping property (5.7) for  $|s|<\frac12$,
 because $\widehat q^+$ is a bounded symbol, 
and $e^+\ol
H^{s}(\rnp)$ identifies with $\dot H^{s}(\crnp)$ then. Now let
 $s=\frac32-\varepsilon $ for a small $\varepsilon >0$. Here we need to show that when 
$u\in \ol H^{\frac32-\varepsilon }(\rnp) $, then $r^+\partial_j\widehat Q^+e^+u\in \ol
H^{\frac12-\varepsilon  }(\rnp)$ for $j=1,\dots,n$. For $j<n$, this follows
simply
because $\partial_j$ can be commuted through $r^+$, $\widehat Q^+$ and $e^+$ so
that we can use that $\partial_ju\in \ol H^{\frac12-\varepsilon }(\rnp)$. For $j=n$, we
proceed as follows:

Since $u\in \ol H^{\frac32-\varepsilon }(\rnp)$,
$$
\partial_ne^+u=e^+\partial_nu+(\gamma _0u)(x')\otimes \delta (x_n),\quad \gamma
_0u\in H^{1-\varepsilon }(\R^{n-1}),\tag5.8
$$
and therefore, since $\widehat Q^+=I+\widehat G^+$ where $\widehat
G^+=\Op(\widehat g^+(\xi ))$ from Proposition 4.3,
$$
\partial_n\widehat Q^+e^+u=\widehat Q^+\partial_ne^+u=\widehat Q^+e^+\partial_nu+(I+\widehat G^+)(\gamma _0u\otimes \delta (x_n)).
$$
In the restriction to $\rnp$, $r^+I(\gamma _0u\otimes \delta (x_n))$ drops
out, so we are left with
$$
r^+\partial^n\widehat Q^+e^+u=r^+\widehat Q^+\partial_ne^+u=r^+\widehat
Q^+e^+\partial_nu+K_{\widehat g^+}\gamma _0u,\quad K_{\widehat g^+}\varphi  =r^+\widehat G^+(\varphi (x')\otimes\delta (x_n)).
$$
Here $K_{\widehat g^+}$ is a potential operator (in the terminology of Eskin
\cite{E81} and Rempel-Schulze \cite{RS82}, generalizing the concept of
Poisson operator of Boutet de Monvel \cite{B66,B71}), which acts as follows:
$$
K_{\widehat g^+}\varphi =r^+\F^{-1}[\widehat g^+(\xi )\hat\varphi (\xi
')].
$$
By (4.8),
$$
|\widehat g^+(\xi )|\le C\ang{\xi '}^{1-\varepsilon /2}\ang{\xi
 }^{\varepsilon /2-1},
 $$
hence
$$
\aligned
\|K_{\widehat g^+}\varphi \|^2_{\ol H^{\frac12-\varepsilon }(\rnp)}&\le \|\widehat G^+(\varphi \otimes \delta )\|^2_{H^{\frac12-\varepsilon }(\rn)}=c\int_{\rn}|\widehat g(\xi )|^2|\hat
\varphi (\xi ')|^2\ang\xi ^{1-2\varepsilon }\,d\xi \\
&\le C
\int_{\rn}|\hat \varphi (\xi ') |^2\ang\xi ^{1-2\varepsilon -2+\varepsilon
}\ang{\xi '}^{2-\varepsilon }\,d\xi 
=C\int_{\rn}|\hat \varphi (\xi ')
|^2\ang\xi ^{-1-\varepsilon  }\ang{\xi '}^{2-\varepsilon
}\,d\xi\\
&
=C'\int_{\R^{n-1}}|\hat \varphi (\xi ') |^2\ang{\xi '}^{2-2\varepsilon
 }\,d\xi '=C''\|\varphi \|^2_{H^{1-\varepsilon }(\R^{n-1}) },
\endaligned
$$
since $\int_{\R}\ang \xi ^{-1-\varepsilon }\,d\xi _n=\ang{\xi
'}^{-\varepsilon }\int_{\R}\ang{\eta _n}^{-1-\varepsilon }\,d\eta _n$.
Inserting $\varphi =\gamma _0u$, we thus have
$$
\|K_{\widehat g^+}\gamma _0u \|_{\ol H^{\frac12-\varepsilon
}(\rnp)}\le C_1\|\gamma _0u \|_{H^{1-\varepsilon }(\R^{n-1}) }\le
C_2\|u\|_{\overline H^{\frac32-\varepsilon }(\rnp)}.
$$
Thus
$$
\|r^+\partial^n\widehat Q^+e^+u\|_{\ol H^{\frac12-\varepsilon
}}\le \|r^+\widehat
Q^+e^+\partial_nu\|_{\ol H^{\frac12-\varepsilon
}}+\|K_{\widehat g^+}\gamma _0u\|_{\ol H^{\frac12-\varepsilon
}}\le C_3\|u\|_{\ol H^{\frac32-\varepsilon
}}.
$$

Altogether, this shows the desired mapping property for $s=\frac32-\varepsilon $, and the
property for general $\frac12\le s <\frac32$ follows by interpolation with the
case $s=0$.

The mapping property (5.7) holds for the inverse $(\widehat Q^+)^{-1}$, since its
symbol $(q^+)^{-1}$ equals $1+\sum_{k\ge 1}(-b_+)^k$ with essentially the same structure.

The identity $((\widehat Q^+)^{-1})_+\widehat Q^+_+=I=\widehat
Q^+_+((\widehat Q^+)^{-1})_+$ valid on $L_2(\rnp)$, holds a fortiori on
$\ol H^s(\rnp)$ for $0< s<\frac32$, and extends by continuity to $\ol H^s(\rnp)$ for $-\frac12< s<0$.
\qed 
\enddemo
 
When $P$ merely satisfies Assumption 3.2, we can still show a useful
forward mapping property of $\widehat Q$, based on the decomposition
in Proposition 4.2.

\proclaim{Proposition 5.2} Let $P$ satisfy Assumption {\rm 3.2}, and
consider $\widehat Q$ and $\widehat F_\pm=\Op(\widehat f_\pm)$ derived from it in Section {\rm 3}.

The operator $\widehat F_{+,+}=r^+\widehat F_+e^+$ is continuous
$$
r^+\widehat F_+e^+\colon \ol H^{s}(\rnp)\to \ol H^{s }(\rnp)\text{ for any }-\tfrac12<s<\tfrac32.\tag5.9
$$

The operator $\widehat F_{-,+}=r^+\widehat F_-e^+$ is continuous from
$ \ol H^{s}(\rnp)$ to $\ol H^{s }(\rnp)$ for
any $s\in \R$.

The operator $\widehat Q_{+}=r^+\widehat Qe^+$ is continuous
$$
r^+\widehat Qe^+\colon \ol H^{s}(\rnp)\to \ol H^{s }(\rnp)\text{ for any }-\tfrac12<s<\tfrac32.\tag5.10
$$
\endproclaim

\demo{Proof} Since $\widehat F_+$ has bounded symbol, it maps $\dot
H^s(\crnp)$ into $H^s(\rn)$ for all $s$, so for $|s|<\frac12$, (5.9)
follows since $\dot
H^s(\crnp)=e^+\ol H^s(\rnp)$ then. For $\frac12<s<\frac32$, we proceed
as in the proof of Proposition 5.1, using that
$$
r^+\partial^n\widehat F_+e^+u=r^+\widehat F_+\partial_ne^+u=r^+\widehat
F_+e^+\partial_nu+K_{\widehat f_+}\gamma _0u,\quad K_{\widehat f_+}\varphi  =r^+\widehat F_+(\varphi (x')\otimes\delta (x_n)),
$$
where $K_{\widehat f_+}$ satisfies similar estimates as $K_{\widehat
g^+}$ by Proposition 4.2.

For $r^+\widehat F_-e^+$, the statement follows since it is on $\rnp$ the
adjoint of $\Op(\overline{\widehat f_-})$, which preserves support in
$\crnp$ and therefore maps $\dot H^s(\crnp)$ into itself for all $s\in
\R$. For $\widehat Q$, the statement now follows since it equals
$s_0+\widehat F_-+\widehat F_+$.\qed
\enddemo

This is as far as we get by applying Lemma 4.1 to $f$. To obtain the
mapping property for higher $s$ would require a control over the
potential operators
$$
\varphi \mapsto r^+\Op(\xi _n^j\widehat f_+(\xi ))(
\varphi (x')\otimes \delta (x_n))
$$
for $j\ge 1$ as well. At any rate, the property shown in Proposition 5.2 will be
sufficient for the integration by parts formulas we are aiming for.

In the elliptic case, we conclude from Proposition 5.1 for the operator $\widehat Q$:

\proclaim{Corollary 5.3} Let $P$ satisfy Assumption {\rm 3.1}, and
consider the operators $\widehat Q,\widehat Q^+, \widehat Q^-$ with
symbols $\widehat q,\widehat q^+, \widehat q^-$ derived from it in Section {\rm 3}.
The operator $\widehat Q_+\equiv r^+\widehat Qe^+$ acts like $r^+\widehat
Q^-e^+r^+\widehat Q^+e^+=\widehat Q^-_+\widehat Q^+_+$, mapping continuously and bijectively
$$
\widehat Q_+=r^+\widehat Qe^+\colon \ol H^{s}(\rnp)\simto \ol H^{s }(\rnp)\text{ for
}-\tfrac12<s< \tfrac32, \tag5.11
$$
and the inverse (continuous in the opposite direction) equals
$$
( r^+\widehat Qe^+)^{-1}=r^+(\widehat
Q^+)^{-1}e^+r^+(\widehat Q^-)^{-1}e^+. \tag5.12
$$

\endproclaim

\demo{Proof} We have for $u\in \dot H^s(\crnp)\simeq \ol H^s(\rnp)$, $|s|<\frac12$, that 
$$
r^+\widehat Qe^+u=r^+\widehat Q^-\widehat Q^+e^+u =r^+\widehat Q^-(e^+r^++e^-r^-)\widehat Q^+e^+u=r^+\widehat Q^-e^+r^+\widehat Q^+e^+u,
$$
 since $r^-\widehat Q^+e^+u=0$; this identity is also valid on the
 subspaces  $\ol H^s(\rnp)$ with $s\ge \frac12$. Combining the
 homeomorphism property of $r^+\widehat Q^+e^+$ shown in Proposition 5.1
 with the known homeomorphism property of  $r^+\widehat Q^-e^+$ on $\ol
 H^s(\rnp)$-spaces (cf.\ (5.6)), we get (5.11). The inverse is pinned down by using
 that $r^+\widehat Q^-e^+$ has inverse $r^+(\widehat Q^-)^{-1}e^+$ on
 $\ol H^s(\rnp)$ for all $s$, and $r^+\widehat Q^+e^+$ has inverse $r^+(\widehat Q^+)^{-1}e^+$ on
 $\ol H^s(\rnp)$ for  $-\frac12<s<\frac32$ in view of Proposition 5.1.\qed 
\enddemo

\subsubhead 5.3 Mapping properties of $\widehat P$ using $\mu $-transmission spaces \endsubsubhead

Now turn the attention to $\widehat P$, which is related to  $\widehat Q$ by
$$
\widehat P=\Xi _-^{\mu '}\widehat Q \,\Xi _+^\mu ,\quad \widehat Q=\Xi _-^{-\mu '}\widehat P \,\Xi _+^{-\mu},\tag5.13
$$
cf.\ (3.12)--(3.13).

We shall describe the solutions of the homogeneous Dirichlet problem
(in the strongly elliptic case)
$$
r^+\widehat Pu=f,\quad \supp u\subset \crnp,\tag5.14
$$
with $f$ given in a space $\ol H^{s }(\rnp)$, and $u$ assumed a priori
to lie in a space $\dot H^\sigma (\comega)$ for low $\sigma $, e.g.\ with $\sigma =a$.

First we observe for  $\Xi ^{-\mu '}_{-,+}=r^+\Xi ^{-\mu
'}_{-}e^+$ that
$$
\Xi ^{-\mu '}_{-,+}r^+\widehat
P=r^+\Xi ^{-\mu '}_{-}\widehat P, 
\tag5.15
$$
since,
as accounted for in \cite{G15, Rem.\
1.1, (1.13)}, the action of $r^+\Xi ^{-\mu '}_{-}$ is independent of
how $r^+\widehat P$ is extended into $\crnm$. Thus, in view of
the mapping properties (5.2) of $\Xi ^{-\mu '}_{-,+}$,
$$
\|r^+\widehat Pu\|_{\ol H^s(\rnp)}\simeq \|\Xi ^{-\mu '}_{-,+}r^+\widehat Pu\|_{\ol H^{s+\operatorname{Re}\mu '}(\rnp)}=\|r^+\Xi _-^{-\mu '}\widehat Pu\|_{\ol H^{s+\operatorname{Re}\mu '}(\rnp)}.\tag5.16
  $$
Composing the equation in (5.14)  with $\Xi _{-,+}^{-\mu '}$ to the left, we can therefore
write it as 
$$
r^+\Xi _-^{-\mu '}Pu=g, \text{ where }g=\Xi _{-,+}^{-\mu '}f\in \ol H^{s+\operatorname{Re}\mu '}(\rnp).\tag5.17
$$

Next, we shall also replace $u$.  Because of the right-hand factor $\Xi
_+^{-\mu }$ in the expression for $\widehat Q$ in (5.13), we need to introduce the $\mu $-transmission spaces
$$
H^{\mu (t)}(\crnp)\equiv \Xi _+^{-\mu }e^+\ol H^{t-\operatorname{Re}\mu
}(\rnp)\text{ for }t>\operatorname{Re}\mu -\tfrac12,\tag 5.18
$$
defined in \cite{G15}; they are Hilbert spaces. (For $t\le
\operatorname{Re}\mu -\tfrac12 $, the convention is to take $H^{\mu
(t)}(\crnp)=\dot H^t(\crnp)$, but this is rarely used.) The following
properties were shown in \cite{G15}:

\proclaim{Theorem 5.4} \cite{G15} Let $t>\operatorname{Re}\mu -\tfrac12$.

$1^\circ$ The mapping $r^+\Xi _+^\mu $ is a homeomorphism of $ H^{\mu
(t)}(\crnp) $ onto $ \ol H^{t-\operatorname{Re}\mu }(\rnp)$ with
inverse $\Xi _+^{-\mu }e^+$.

$2^\circ$ For $|t-\operatorname{Re}\mu |<\frac12$, $ H^{\mu
(t)}(\crnp) =
\dot H^{t }(\crnp)$.

$3^\circ$ Assume $\operatorname{Re}\mu >-1$ and $t>\operatorname{Re}\mu
+\frac12$. Then
$$
H^{\mu (t) }(\crnp)
\subset \dot H^t(\crnp)+x_n^\mu e^+\ol H^{t-\operatorname{Re}\mu }(\rnp),\tag5.19
$$
where $\dot H^t(\crnp)$ is replaced by $\dot H^{t-\varepsilon
}(\crnp)$ if $t-\operatorname{Re}\mu -\frac12\in \N$. Moreover, the trace of $u/x_n^{\mu }$ is
well-defined on $H^{\mu (t) }(\crnp)$ and satisfies
$$
\Gamma(1+\mu )\gamma
_0(u/x_n^{\mu })=\gamma _0\Xi _+^{\mu }u\in H^{t-\operatorname{Re}\mu -\frac12}(\R^{n-1}).
\tag5.20
$$
\endproclaim

Rule $1^\circ$ is shown in \cite{G15, Prop.\ 1.7}. 
 Rule $2^\circ$, shown in \cite{G15,
(1.26)}, holds because of  the mapping property (5.2) for $\Xi
 _+^{-\mu }$ and the identification of $ e^+\ol
H^{t-\operatorname{Re}\mu }(\rnp)$ with $\dot
H^{t-\operatorname{Re}\mu }(\crnp)$ when $t-\operatorname{Re}\mu\in\,]-\frac12,\frac12[\,$.
Rule $3^\circ$ is shown in \cite{G15, Th.\ 5.1, Cor.\ 5.3,
Th.\ 5.4}; it deals with a higher $t$, where $ e^+\ol
H^{t-\operatorname{Re}\mu  }(\rnp)$ has a jump at $x_n=0$, and the
 coefficient $x_n^\mu $ appears. Let us just mention the key formula
$$
{\Cal F}^{-1}_{\xi _n\to x_n}[(\ang{\xi '}+i\xi _n)^{-\mu }(\ang{\xi '}+i\xi _n)^{-1}]
=\tfrac
1{\Gamma (\mu +1)}e^+r^+x_n^{\mu }e^{-\ang{\xi '} x_n},
$$
which indicates how $\Xi _+^{-\mu }=\Op((\ang{\xi '}+i\xi _n)^{-\mu })$ is connected
with the factor $x_n^\mu $. Besides in \cite{G15, Sect.\ 5}, explicit calculations are carried out e.g.\ in
\cite{G19, Lemma 3.3} (and \cite{G14, Appendix}).

We note in passing that in the definition (5.18),  one can equivalently replace the order-reducing
operator family $\Xi
_+^{t}=\Op((\ang{\xi '}+i\xi _n)^t)$ by $\Op(([\xi
']+i\xi _n)^t)$, or by $\Lambda _+^t$, as defined in \cite{G15}.

Now continue the discussion of (5.17):
In view of Theorem 5.4 $1^\circ$, we can set $ v=r^+\Xi
_+^\mu u$, where $u=\Xi _+^{-\mu }e^+v$, and hereby
$$
r^+\Xi ^{-\mu '}\widehat Pu=r^+\Xi ^{-\mu '}\widehat P\Xi
_+^{-\mu }e^+v=r^+\widehat Qe^+v =\widehat Q_+v.
$$
Then 
the equation (5.17) reduces to an equivalent equation
$$
\widehat Q_+v=g,
$$
with $g$ given in $ \ol H^{s+\operatorname{Re}\mu '}(\rnp)$ and $v$ a
priori taken in $ \dot H^{\sigma -\operatorname{Re}\mu }(\crnp)$.
This was solved in Corollary 5.3, so we find for $r^+\widehat P$:

\proclaim{Theorem 5.5} Let $P$ satisfy Assumption {\rm 3.1}.
For $\operatorname{Re}\mu -\frac12<t< \operatorname{Re}\mu +\frac32$,
$r^+\widehat P$ defines a homeomorphism (continuous bijective operator
with continuous inverse)
$$
r^+\widehat P\colon  H^{\mu (t) }(\crnp)\simto \ol H^{t -2a }(\rnp). \tag5.21
$$

Furthermore, if $u$ is in $ \dot H^\sigma (\crnp)$ for some $\sigma >\operatorname{Re}\mu
-\frac12$ (this includes the value $\sigma =a$) and solves {\rm (5.14)},
then $u\in  H^{\mu (t) }(\crnp)$.

Here $\operatorname{Re}\mu >-\frac12$ since $a>0$ and $|\operatorname{Re}\delta
|<\frac12$, so the rules in {\rm
Theorem 5.4 $3^\circ$} apply.

\endproclaim

\demo{Proof} In view of (5.15) and (5.16), and the mapping property
of $\widehat Q$ established in Corollary 5.3, $r^+\widehat P$ has the forward
mapping property in (5.21).

To solve (5.14), let $\sigma =\operatorname{Re}\mu
-\frac12+\varepsilon $ for a small $\varepsilon $,  set $g=\Xi
_{-,+}^{-\mu '}f\in \ol H^{t-2a+\operatorname{Re}\mu '}(\rnp)=\ol
H^{t-\operatorname{Re}\mu }(\rnp)$ and $v=r^+\Xi _+^\mu u\in \dot
H^{-\frac12+\varepsilon  }(\crnp)$. Then (5.14)
reduces to solving 
$$
\widehat Q_+v=g,\tag5.22
$$
with $g$ given in $\ol H^{t-\operatorname{Re}\mu }(\rnp)$ and
$v$ a priori lying in $\dot H^{-\frac12+\varepsilon  }(\crnp)$.
By Corollary 5.3, (5.22) has a unique solution  $v\in \ol
H^{t-\operatorname{Re}\mu }(\rnp)$, so $u$ must lie in  $ H^{\mu
(t)}(\crnp)$, and the mapping (5.21) is bijective.
 \qed
\enddemo

\example{Remark 5.6} This theorem differs from the strategy
pursued in \cite{E81}, and gives a new insight. The technique in \cite{E81, Th.\ 7.3} for showing solvability in
a higher-order Sobolev space, say with $\frac12<t-\operatorname{Re}\mu <\frac32$, $f$ given in $\ol
H^{t-2a}(\rnp)$, is 
to supplement $\widehat P$ with a potential operator $K_{\widehat P}$ constructed
from $\widehat P$ such that the solutions are of the form
$u=u_++K_{\widehat P}\varphi $
with $u_+\in \dot H^t(\crnp)$, $\varphi $ a generalized trace  derived from $f$. Our aim is
to show that there is a {\it universal
description} of the space of solutions $u$ of (5.14)  with right-hand side in $ \ol H^{t-2a}(\rnp)$, that
depends only on $\mu $, and applies to any $P$ of the given type. The
$\mu $-transmission spaces (5.18) serve this purpose. In \cite{G15}, they are
shown to have  
this role for arbitrarily high $t$  when the full $\mu
$-transmission condition holds.
\endexample

One more important property of $\mu $-transmision spaces is that the
 spaces with $C^\infty $-ingredients
 $\E_\mu (\crnp)\cap \Cal E'(\rn)$ and $e^+x_n^\mu \Cal S(\crnp)$ {\it
 are
 dense subsets of} $H^{\mu
(t)}(\crnp)$ for all $t>\operatorname{Re}\mu -\frac12$,
 $\operatorname{Re}\mu >-1$ (cf.\ \cite{G15, Prop.\ 4.1} and \cite{G21, Lemma 7.1}).
Recall also (5.5), which makes the statement for $e^+x_n^\mu \Cal
 S(\crnp)$ rather
 evident, since $\Cal S(\crnp)$ is dense in $\ol H^s(\rnp)$ for all $s\in\R$.
 Hence $r^+\widehat P$ applies
nicely to these spaces.

When $P$ merely satisfies Assumption 3.2, we have at least the forward
mapping part of (5.21):

\proclaim{Theorem 5.7}  Let $P$ satisfy Assumption {\rm 3.2}.
For $\operatorname{Re}\mu -\frac12<t< \operatorname{Re}\mu +\frac32$,
$r^+\widehat P$ maps continuously
$$
r^+\widehat P\colon  H^{\mu (t) }(\crnp)\to \ol H^{t -2a }(\rnp). \tag5.23
$$
\endproclaim

\demo{Proof} This follows as in the preceding proof, now using the mapping property of
$r^+\widehat Qe^+$ established in Proposition 5.2.\qed
\enddemo

\subsubhead 5.4 Consequences for the given operator $P$ \endsubsubhead

The following consequences can be drawn for the original operator $P$:

\proclaim{Theorem 5.8} 

$1^\circ$ Let $P$ satisfy Assumption {\rm 3.2}.  Then $P=\widehat
P+P'$, where $\widehat P$ is defined by {\rm (3.9)} and $P'$ is of
order $2a-1$. For $\operatorname{Re}\mu -\frac12<t<\operatorname{Re}\mu +\frac32$, $r^+P$ maps continuously
$$
r^+ P\colon  H^{\mu (t) }(\crnp)\to \ol H^{t -2a }(\rnp). \tag5.24
$$

$2^\circ$ Let $P$ satisfy Assumption {\rm 3.1}. Then in the
decomposition  $P=\widehat
P+P'$, $r^+\widehat P$ is invertible, as described in Theorem {\rm
5.5}.

Let $\operatorname{Re}\mu -\frac12<t<\operatorname{Re}\mu
+\frac32$, let $f\in
\ol H^{t-2a}(\rnp)$, and let  $u\in \dot H^{\sigma  }(\crnp) $ (for
some $\sigma >\operatorname{Re}\mu -\frac12$) solve the homogeneous
Dirichlet problem
$$
r^+Pu=f \text{ on }\rnp, \quad \supp u\subset \crnp.\tag5.25
$$
Then $u\in H^{\mu (t)}(\crnp)$.
\endproclaim

\demo{Proof} The original operator $P$ equals $\operatorname{Op}(p(\xi
))$ with $p(\xi )$  homogeneous on $\rn$ of degree $2a>0$; in
particular it is continuous at 0. It is decomposed into 
$$
p(\xi )=\widehat p(\xi )+p'(\xi ).\tag5.26
$$
where  $p'(\xi )$ is $O(\ang\xi ^{2a-1})$ for $|\xi |\ge 2$ by (3.10)
and continuous, hence
$$
|p'(\xi )|\le C'\ang\xi ^{2a-1}\text{ for }\xi \in \rn.
$$
This implies that $P'=\Op(p')$ maps $H^s(\rn)$ continuously into
$H^{s-2a+1}(\rn)$ for all $s\in\R$, and hence
$$
r^+P'\colon \dot H^s(\crnp)\to
\ol H^{s-2a+1}(\rnp)\text{ for all }s\in\R.\tag5.27 
$$

$1^\circ$. The forward mapping property (5.23) holds for $r^+\widehat P$ by
Theorem 5.7. To show that it holds for
$r^+P'$, let $\operatorname{Re}\mu -\frac12<t< \operatorname{Re}\mu +\frac32$. 

If $t-\operatorname{Re}\mu<\frac12$, then $H^{\mu (t) }(\crnp)=\dot
H^t(\crnp)$, and $r^+P'\dot H^t(\crnp)\subset \ol
H^{t-2a+1}(\rnp)\subset \ol H^{t-2a}(\rnp)$ by (5.27), matching the mapping
property of $\widehat P$.

If $\frac12\le t-\operatorname{Re}\mu < \frac32$, we use the definition of
$H^{\mu (t)}(\crnp)$ to see that for small  $\varepsilon >0$,
$$
\aligned
r^+P' H^{\mu (t)}(\crnp)&=r^+P'\Xi _+^{-\mu }e^+\ol H^{t-  \operatorname{Re}\mu }(\rnp)
\subset r^+P'\Xi _+^{-\mu }\dot H^{\frac12-\varepsilon  }(\crnp)\\
&=r^+P'\dot H^{\frac12-\varepsilon +\operatorname{Re}\mu }(\crnp)
\subset \ol H^{\frac32-\varepsilon +\operatorname{Re}\mu -2a
}(\crnp)\subset \ol H^{t-2a}(\rnp),
\endaligned
$$
also matching the mapping property of $\widehat P$.

Now (5.24)
follows by adding the statements for $P'$ and $\widehat P$. This
shows $1^\circ$.

$2^\circ$. The first statement registers what we already know about
$r^+\widehat P$. Proof of the regularity statement: With $u$ and $f$ as
defined there, denote  $\sigma =\operatorname{Re}\mu
-\frac12+\varepsilon  $; here $\varepsilon >0$.
 Then 
$$
r^+\widehat Pu=r^+Pu-r^+P'u\in \ol H^{t-2a}(\rnp)+\ol H^{\operatorname{Re}\mu
+\frac12+\varepsilon-2a}(\rnp).
$$

If $t\le \operatorname{Re}\mu +\frac12 +\varepsilon  $, $r^+\widehat Pu\in  \ol
H^{t-2a}(\rnp) $, and we conclude from Theorem 5.5 that $u\in H^{\mu
(t)}(\crnp)$.

If $t> \operatorname{Re}\mu +\frac12 +\varepsilon  $, $r^+\widehat Pu\in  \ol
H^{\operatorname{Re}\mu +\frac12+\varepsilon -2a}(\rnp) $; here
Theorem 5.5 applies to  give the
intermediate information that  $u\in H^{\mu
(\operatorname{Re}\mu +\frac12+\varepsilon )}(\crnp)$. From this follows that
$$
u\in \Xi _+^{-\mu }e^+\ol H ^{\frac12+\varepsilon }(\rnp)\subset
\Xi _+^{-\mu }\dot H ^{\frac12-\varepsilon '}(\crnp)=\dot H ^{\operatorname{Re}\mu +\frac12-\varepsilon '}(\crnp),
$$
for any $\varepsilon '>0$. Then $r^+P'u\in \ol H ^{\operatorname{Re}\mu
+\frac32-\varepsilon '-2a}(\rnp) $. Choosing $\varepsilon '$ so small
that  $\operatorname{Re}\mu
+\frac32-\varepsilon '\ge t$, we have that  $r^+P'u\in \ol H
^{t-2a}(\crnp) $; hence  $r^+\widehat Pu\in \ol H
^{t-2a}(\crnp) $, so it follows from Theorem 5.5 that $u\in  H
^{\mu (t)}(\crnp) $. This ends the proof of $2^\circ$.\qed
\enddemo

\example{Example 5.9} Theorem 5.8 applies to the operator
$L=\operatorname{Op}(\Cal L(\xi ))$ described in (3.5)ff.,  showing
that it maps $ H^{\mu (t) }(\crnp)$ to $\ol H^{t -2a }(\rnp)$ for
$-\frac12<t-\mu <\frac32$, and that solutions of the homogeneous
Dirichlet problem with $f\in \ol H^{t -2a
}(\rnp)$ are in $ H^{\mu (t) }(\crnp)$ for these $t$. The appearance
of the factor $x_n^\mu $ (cf.\ (5.19)) is consistent with the
regularity shown in terms of H\"older spaces in \cite{DRSV21}.

In particular, the result
provides a valid basis for applying $r^+L$ to $\E_\mu (\crnp)\cap \Cal
E'(\rn)$ or $e^+x_n^{\mu }\Cal S(\crnp)$, mapping these spaces into $\bigcap_{\varepsilon >0}\ol
H^{\frac32-a+\delta -\varepsilon }(\rnp)$.
\endexample

\example{Remark 5.10} The domain spaces $ H^{\mu (t) }(\crnp)$ entering in
Theorem 5.8 can be precisely described: For
$|t-\operatorname{Re}\mu |<\frac12$, we already know from Theorem 5.4 $2^\circ$ that
$ H^{\mu (t) }(\crnp)=\dot H^t(\crnp)$. For
$\frac12<t-\operatorname{Re}\mu <\frac32$, we have by \cite{G19, Lemma
3.3} that $u\in  H^{\mu (t) }(\crnp)$ if and only if 
$$
u=v+w,\text{ where }w\in \dot H^t(\crnp)\text{ and }v=e^+x_n^\mu K_0\gamma _0(u/x_n^\mu );
$$
here $K_0$ is the Poisson operator $K_0\colon \varphi \mapsto z$ solving the Dirichlet problem for $1-\Delta $,
$$
(1-\Delta )z=0\text{ on }\rnp,\quad \gamma _0 z=\varphi \text{ at }x_n=0,
$$
with $\varphi \in H^{t-\operatorname{Re}\mu -\frac12}(\R^{n-1})$. For
$t-\operatorname{Re}\mu =\frac12$, we have the information $u\in \bigcap_{\varepsilon
>0}\dot H^{t-\varepsilon }(\crnp)$. As a concrete example, the elements
$u$ of $ H^{\frac12(\frac32)}(\crnp)$ are the functions $u=v+w$, where
$w\in \dot H^\frac32(\crnp)$ and $v=x_n^\frac12 K_0\varphi $ for some $\varphi
\in H^\frac12(\R^{n-1})$; this $\varphi $ equals $\gamma _0(u/x_n^\frac12)$.

\endexample

\subhead 6. The integration by parts formula \endsubhead

\subsubhead 6.1 An integration by parts formula for $\widehat P$ \endsubsubhead

It will now be shown that the operators $P$ satisfying merely the
principal $\mu $-transmission condition (Assumption 3.2) have an
integration by parts formula over $\rnp$, involving traces $\gamma
_0(u/x_n^{\mu })$.
The study will cover the special
operator $L$ in Example 5.9 (regardless of whether a full $\mu
$-transmission condition might hold, as assumed in \cite{G21}). It
also covers more general strongly elliptic operators, and it covers operators
that are not necessarily elliptic.

The basic observation is:

\proclaim{Proposition 6.1} Let $\mu  \in\C$.
Let $w\in \bigcap_{s}\ol H^s(\rnp)$, 
and let
$u'\in \E_{\bar\mu }(\crnp)\cap
\E'(\rn)$. Denote  $w'=r^+\Xi _+^{\bar\mu }u'\in  \bigcap_{s}\ol H^s(\rnp) $; correspondingly $u'=\Xi
_+^ {-\bar\mu }e^+w'$ in view of {\rm Theorem 5.4 $1^\circ$}. Then
$$
(I\equiv )\int _{\rnp}\Xi _-^{ \mu }e^+w\,\partial_n\bar u'\, dx=
(\gamma _0w,\gamma _0w')_{L_2(\R^{n-1})}+(w,\partial_nw')_{L_2(\rnp)}
.\tag 6.1
$$
The left-hand side is interpreted as in {\rm (6.2)} below when
$\operatorname{Re}\mu \le 0$.

The formula extends to $w\in \ol H^{\frac12+\varepsilon}(\rnp)$  and
$u'\in H^{\bar\mu (t)}(\crnp)$ with $t\ge\operatorname{Re}\mu
+\frac12-\varepsilon $ (any small $\varepsilon >0$),
using the representation {\rm (6.2)}.
\endproclaim

\demo{Proof}
This was proved in \cite{G16, Th.\ 3.1} for $\mu =a>0$ (see also
Remark 3.2 there with the elementary case $a=1$), and in  \cite{G21,
Th.\ 4.1} for real $\mu >-\frac12$, so the main task is to check
that the larger range of complex $\mu $ is allowed. We write $\overline{u'}$ as $\bar u'$ for short.

Note that  when $\operatorname{Im}\mu \ne 0$,  $\E_{\mu }(\crnp)$ is
different from $\E_{\operatorname{Re}\mu }(\crnp)$, e.g.\  since $x_n^\mu
/x_n^{\operatorname{Re}\mu }=x_n^{i\operatorname{Im}\mu
}=e^{i\operatorname{Im}\mu \log x_n}$ has absolute value 1 and is
$C^\infty $ for $x_n>0$, but oscillates when $x_n\to 0$.

By the mapping properties of $\Xi _{-,+}^\mu $ (cf.\ (5.2)),  $r^+\Xi _{-}^\mu e^+w\in
\bigcap_s\ol H^s(\rnp) $, hence is integrable.   
When $\operatorname{Re}\mu >0$, the function $\partial_nu'$ is
$O(x_n^{\operatorname{Re}\mu -1})$ and compactly supported, so 
 the
left-hand side of (6.1) makes sense as an integral of an
$L_1$-function. When $\mu $ is general, we observe 
that for any small $\varepsilon >0$,
$$\partial_nu'\in \E_{\bar\mu -1}(\crnp)\cap \E'(\rn)\subset \Xi
_+^{1-\bar\mu }e^+\bigcap_s\ol H^s(\rnp)\subset \Xi
_+^{1-\bar\mu }\dot H^{\frac12-\varepsilon }(\crnp)=\dot H^{-\frac12+\operatorname{Re}\mu -\varepsilon }(\crnp),
$$
so the integral $I$ makes sense as the duality
$$
I=\ang{r^+\Xi _-^{ \mu }e^+w,\partial_n u'}_{\ol H^{\frac12-\operatorname{Re}\mu +\varepsilon }(\rnp),\dot H^{-\frac12+\operatorname{Re}\mu -\varepsilon }(\crnp)}.\tag6.2
$$
Since the adjoint of $r^+\Xi _-^{ \mu }e^+$ equals $\Xi _+^{\bar \mu
}$, $I$ is by transposition turned into
$$
I=\ang{w,\Xi _+^{\bar\mu }\partial_n u'}_{\ol H^{\frac12 +\varepsilon
},\dot H^{-\frac12 -\varepsilon }}
=\ang{w,\partial_n\Xi _+^{\bar\mu }\Xi _+^{-\bar\mu }e^+w'}_{\ol H^{\frac12 +\varepsilon
},\dot H^{-\frac12 -\varepsilon }}
=\ang{w,\partial_ne^+w'}_{\ol H^{\frac12 +\varepsilon },\dot H^{-\frac12 -\varepsilon }}.
$$
Note that $\partial_ne^+w$ satisfies an equation like (5.8), which fits
in here since the space $\dot H^{-\frac12 -\varepsilon }(\crnp)$ contains
distributions of the form $\varphi (x')\otimes \delta (x_n)$. The expression is analysed as in \cite{G16 Th. 3.1} (and \cite{G21 Th.\ 4.1}), leading to
$$
I=(\gamma _0w,\gamma _0 w')_{L_2(\R^{n-1})}+(w,e^+\partial_nw')_{L_2(\rnp)},\tag6.3
$$
which shows (6.1).

For the whole analysis, it suffices that $w\in \ol H^s(\rnp)$ with
$s=\frac12+\varepsilon $, since $\Xi ^\mu _{-,+}w\in \ol
H^{\frac12-\operatorname{Re}\mu +\varepsilon }(\rnp)$ then. For $u'$, it
then suffices that $u'\in H^{\bar\mu (t)}(\crnp)=\Xi _+^{-\bar\mu }e^+
\ol H^{t-\operatorname{Re}\mu }(\rnp)$ with $t\ge \operatorname{Re}\mu
+\frac12-\varepsilon $ (assuming $0<\varepsilon <1$), since 
$$\partial_nu'\in  \Xi
_+^{1-\bar\mu }e^+\ol H^{t-\operatorname{Re}\mu }(\rnp)\subset \Xi
_+^{1-\bar\mu }\dot H^{\frac12-\varepsilon   }(\crnp)
=\dot H^{-\frac12+\operatorname{Re}\mu -\varepsilon }(\crnp)
$$
then, so that the duality in (6.2) is well-defined.\qed

\enddemo

We shall now show:

\proclaim{Theorem 6.2} Let $P$ satisfy Assumption {\rm 3.2}; it is of
 order $2a$ and satisfies the principal $\mu $-transmission condition
 in the direction $(0,1)$ for
some $\mu =a+\delta \in\C$, and we denote $a-\delta =\mu '$. Assume moreover that $\operatorname{Re}\mu>-1 $, $\operatorname{Re}\mu'>-1 $.
 Consider  $\widehat P=\Op(\widehat p(\xi ))$,
 as defined by {\rm (3.9)}. 
 For $u\in \E_\mu (\crnp)\cap \E'(\rn)$, $u'\in \E_{\bar\mu '} (\crnp)\cap \E'(\rn)$,
 there holds 
$$
\aligned 
\int_{\rnp} \widehat Pu\,\partial_n\bar
u'\,dx&+\int_{\rnp}\partial_nu\,\overline{  \widehat P^*u'}\,dx\\
&=\Gamma (\mu +1){\Gamma(\mu '+1)}\int_{{\Bbb
R}^{n-1}}s_0\gamma _0(u/x_n^{\mu })\,{\gamma _0(\bar u'/x_n^{\mu '})}\,
dx',
\endaligned
\tag6.4
$$
where $s_0=e^{-i\pi \delta }p(0,1)$. The formula extends to $u\in H^{\mu (t)}(\crnp)$,  $u'\in H^{\bar\mu
'(t')}(\crnp)$, for  $t>\operatorname{Re}\mu +\frac12$,
$t'>\operatorname{Re}\mu '+\frac12$.

The integrals over $\rnp$ are interpreted as dualities (as in Proposition {\rm 6.1}) when $\operatorname{Re}\mu $ or
$\operatorname{Re}\mu '\le 0$, and when extended to general $u,u'$.
\endproclaim

\demo{Proof} Since integration over $\rnp$ in itself indicates that the
functions behind the integration sign are restricted to $\rnp$, we can
leave out the explicit mention of $r^+$.
Recall that
$$
\widehat p= \chi _-^{\mu '}\widehat q\, \chi _+^\mu ,\quad \widehat P=\Xi _-^{\mu '}\widehat Q\,\Xi _+^\mu ,
$$
cf.\ (3.13). The adjoint is $\widehat P^*=\Xi _-^{\bar\mu }\widehat
Q\,\Xi _+^{\bar\mu '}$. Recall from Proposition 4.2 that
$$
q(\xi )=s_0+f_+(\xi )+f_-(\xi ),\text{ hence }\widehat Q =s_0+  \widehat F_++ \widehat F_-,
$$
where $\widehat f_\pm(\xi )$ extend holomorphically in $\xi _n+i\tau $
 into $\C_-$ resp.\
 $\C_+$, estimated as in (4.8).

Accordingly, $\widehat P$ splits up in three terms
$$
\widehat P  =\widehat P_1+\widehat P_2+\widehat P_3,\text{ where }
\widehat P_1 =s_0 \Xi _-^{\mu '}\Xi _+^\mu , \quad
\widehat P_2  =\Xi _-^{\mu '}\widehat F_+\Xi _+^\mu,\quad \widehat P_3 =\Xi _-^{\mu '}\widehat F_-\Xi _+^\mu .\tag6.5
$$

Consider the contribution from $\widehat P_1$:
$$
\int _{\rnp}\widehat P_1u\,{\partial_n \bar u'}\, dx+\int
_{\rnp}\partial_nu\,\overline
{\widehat P_1^* u'}\, dx  =s_0\int _{\rnp}\Xi _-^{\mu '}\Xi _+^\mu u\,{\partial_n \bar u'}\, dx+s_0\int
_{\rnp}\partial_nu\,\overline
{\Xi _-^{\bar \mu }\Xi _+^{\bar\mu '} u'}\, dx .
$$
Recall that $s_0=q(0,1)=e^{-i\pi \delta }p(0,1)$ by (3.16); this constant is left
out of the next calculations.

When $u\in \E_\mu (\crnp)\cap \E'(\rn)$, then  $w=r^+\Xi _+^{\mu }u\in \bigcap_s\ol
H^s (\rnp) $. Similarly as in (5.15), $r^+\Xi _-^{\mu '}\Xi _+^{\mu
}u=r^+\Xi _-^{\mu '}e^+r^+\Xi _+^{\mu}u$, which equals $r^+\Xi _-^{\mu '}e^+w$, 
hence lies in 
$\bigcap_s\ol
H^s (\rnp) $ by (5.2). An application of Proposition 6.1 with $ \mu $ replaced
by $\mu '$ gives:
$$
\int _{\rnp}\Xi _-^{\mu '}\Xi _+^\mu u\,{\partial_n \bar u'}\, dx=\int _{\rnp}r^+\Xi _-^{\mu '}e^+w\,\partial_n\bar u'\, dx=(\gamma _0w,\gamma _0w')_{L_2(\R^{n-1})}+(w,\partial_nw')_{L_2(\rnp)},
$$
where $w'=r^+\Xi _+^{\bar \mu '}u'$.

We can apply
the analogous argument to show that the conjugate of $\int
_{\rnp}\partial_nu\,\overline{\Xi _-^{\bar \mu }\Xi _+^{\bar\mu '}u'}\,dx$ satisfies
$$
\int _{\rnp}\Xi _-^{\bar \mu }\Xi _+^{\bar\mu '}u'\,\partial_n\bar u\, dx=(\gamma _0w',\gamma _0w)_{L_2(\R^{n-1})}+(w',\partial_nw)_{L_2(\rnp)};
$$
here $w'=r^+\Xi _+^{\bar\mu '}u'$ and $w=r^+\Xi _+^{\mu }u$ are the same as
the functions defined in the treatment of the first integral.

It follows by addition that
$$
\multline
\int _{\rnp}\Xi _-^{\mu '}\Xi _+^\mu u\,\partial_n\bar u'\, dx+\int
_{\rnp}\partial_nu\,
\overline {\Xi _-^{\bar \mu }\Xi _+^{\bar\mu '} u'}\, dx\\
=2(\gamma _0w,\gamma _0w')_{L_2(\R^{n-1})}+(w,\partial_n 
w')_{L_2(\rnp)}+(\partial_nw, 
w')_{L_2(\rnp)}=(\gamma _0w,\gamma _0w')_{L_2(\R^{n-1})};
\endmultline
$$
in the last step we used that $\int_{\rnp}(w\partial_n\bar w'+\partial_nw \bar
w')\,dx=-\int_{{\Bbb R}^{n-1}}\gamma _0w\gamma _0\bar w'\, dx'$.
Insertion of $\gamma _0w=\gamma _0\Xi _+^{\mu }u=\Gamma(1+\mu )\gamma
_0(u/x_n^{\mu })$ (cf.\ (5.20)),
and similarly  $\gamma _0w'=\gamma _0\Xi _+^{\bar\mu ' }u'
 =\Gamma(1+\bar\mu ')\gamma
_0(u'/x_n^{\bar\mu ' })$, gives (6.4) with $\widehat P$ replaced by
$\widehat P_1$ (using also that $\overline{
\Gamma(1+\bar\mu ')}=\Gamma(1+\mu ')$).

As for extension of the formula to larger spaces, we note that by Proposition 6.1, the calculations for the first integral allow $w\in \ol
H^{\frac12+\varepsilon }(\rnp)$, corresponding to $u\in H^{\mu
(t)}(\crnp)$ with $t=\operatorname{Re}\mu +\frac12+\varepsilon $, and
$u'\in H^{\bar\mu '(t')}(\crnp)$ with
$t'\ge \operatorname{Re}\mu '+\frac12-\varepsilon $. With the
analogous conditions for the calculations of the second integral, we
find altogether that $t>\operatorname{Re}\mu +\frac12$,
$t'>\operatorname{Re}\mu '+\frac12$,
is allowed.

The contributions from $\widehat P_2$ and $\widehat P_3$ will be
treated by variants of this proof, where we show that their boundary
integrals give zero.

Consider $\widehat P_2$. As
in (5.15), we have:
$$
\aligned
r^+\widehat P_2u&=r^+ \Xi _-^{\mu '}\widehat F_+\Xi _+^{\mu }u= r^+\Xi
_-^{\mu '}e^+r^+(\widehat F_+\Xi
_+^{\mu }u),\\
r^+\widehat P_2^*u'&=
r^+ \Xi _-^{\bar\mu }\widehat F_+^*\Xi _+^{\bar\mu '}u'
=r^+( \Xi _-^{\bar\mu }\widehat F_+^*)e^+r^+\Xi _+^{\bar\mu '}u',
\endaligned
$$
where $\widehat F_+^*=\Op(\overline {\widehat f_+})$.
Set 
$$
w=r^+\Xi _+^{\mu }u,\quad w_1=r^+\widehat F_+\Xi
_+^{\mu }u, 
\quad w'=r^+\Xi _+^{\bar\mu '}u' .\tag 6.6
$$
Here when $u\in H^{\mu
(t)}(\crnp)$, $w\in \ol H^{t-\operatorname{Re}\mu }(\rnp) $, and when
$u'\in H^{\bar\mu '
(t')}(\crnp)$, $w'\in \ol H^{t'-\operatorname{Re}\mu '}(\rnp)
$. For $w_1$ we have since  $u=\Xi _+^{-\mu }e^+w$ 
(by Theorem 5.4 $1^\circ$), that
$$
w_1=r^+\widehat F_+\Xi _+^{\mu }u=r^+\widehat F_+\Xi _+^{\mu }\Xi
_+^{-\mu }e^+w=r^+\widehat F_+e^+w \in  \ol H^{t-\operatorname{Re}\mu }(\rnp),
$$
when $-\frac12<t-\operatorname{Re}\mu <\frac32$, by the mapping
property for $\widehat F_+$ established in Proposition 5.2.

We can then apply Proposition 6.1 to the first integral for $\widehat P_2$,
with $\mu $ replaced by $\mu '$, giving when $t-\operatorname{Re}\mu >\frac12$:
$$
\multline
\int_{\rnp}\widehat P_2u\partial_n\bar u'\,dx=\int_{\rnp}\Xi _-^{\mu '}\widehat F_+\Xi ^\mu _+u\partial_n\bar u'\,dx=
\int_{\rnp}\Xi _-^{\mu '}e^+w_1\partial_n\bar u'\,dx\\
=(\gamma _0w_1,\gamma
_0w')_{L_2(\R^{n-1})}+(w_1,\partial_nw')_{L_2(\rnp)}.
\endmultline
\tag6.7
$$
There is a general formula for the trace, entering in Vishik and
Eskin's calculus as well as that of Boutet de Monvel,
$$
\gamma _0v=(2\pi )^{-n}\int_{{\Bbb
R}^{n-1}}e^{ix'\cdot\xi '}\int_{\R}\F(e^+v)\,d\xi_n
d\xi ',
$$
where the integral over $\R$ is read either as an ordinary integral
or, if necessary, as the integral $\int^+$ defined e.g.\ in \cite{G09,
(10.85)} (also recalled in \cite{G18, (A.15),
(A.1)}). Applying this to $w_1$, we find:
$$
\gamma _0w_1=\gamma _0(\widehat F_{+,+}w)=(2\pi )^{-n}\int_{{\Bbb
R}^{n-1}}e^{ix'\cdot\xi '}\int_{\R}\widehat f_+(\xi ',\xi _n)\F(e^+w)\,d\xi
_nd\xi '.\tag6.8
$$
This integral gives 0 for the following reason: Take $w$ in the dense
subspace of $\ol H^{t-\operatorname{Re}\mu }(\rnp)$ of compactly
supported functions in $C^\infty (\crnp)$.
Both $\widehat f_+$ and $\F(e^+w)$ are holomorphic in $\C_-$ as
functions of $\xi _n$, $f_+$ being $O(\ang{\xi _n}^{-1+\varepsilon })$
and $\F(e^+w)$ being $O(\ang{\xi _n}^{-1})$ on $\overline {\C}_-$,
whereby the integrand is $O(\ang{\xi _n}^{-2+\varepsilon })$ there
(and is in $L_1$ on $\R$); then the integral over $\R$ can
be transformed to a closed contour in $\C_-$ and gives 0.

We can then conclude:
$$
\int_{\rnp}\widehat P_2u\partial_n\bar u'\,dx=\int_{\rnp}\Xi _-^{\mu '}\widehat F_+\Xi ^\mu _+u\partial_n\bar u'\,dx=
(w_1,\partial_nw')_{L_2(\rnp)}.
\tag6.9
$$

The other contribution from $\widehat P_2$ is, in conjugated form,
$$
\aligned
\int_{\rnp}\widehat P_2^*u'\partial_n\bar u\,dx&=\int_{\rnp}\Xi _-^{\bar \mu }\widehat F_+^*e^+r^+\Xi _+^{\bar\mu
'}u'\,\partial_n\bar u\, dx\\
&=\ang {r^+\Xi _-^{\bar \mu }\widehat F_+^*e^+r^+\Xi _+^{\bar\mu
'}u',\partial_nu}_{\ol H^{\frac12-\operatorname{Re}\mu +\varepsilon
},\dot H^{-\frac12+\operatorname{Re}\mu -\varepsilon }}\\
&=\ang {r^+\Xi _+^{\bar\mu
'}u',\widehat F_+\Xi _+^{\mu }\partial_nu}_{\ol H^{\frac12 +\varepsilon },\dot H^{-\frac12 -\varepsilon }}
=\ang {r^+\Xi _+^{\bar\mu
'}u',\partial_n\widehat F_+\Xi _+^{\mu }u}_{\ol H^{\frac12 +\varepsilon },\dot H^{-\frac12 -\varepsilon }}\\
&=\ang {r^+\Xi _+^{\bar\mu
'}u',\partial_n\widehat F_+e^+w}_{\ol H^{\frac12
+\varepsilon },\dot H^{-\frac12 -\varepsilon }}
=
\ang {w',\partial_ne^+w_1}_{\ol H^{\frac12
+\varepsilon },\dot H^{-\frac12 -\varepsilon }}\\
&=(\gamma _0w',\gamma
_0w_1)_{L_2(\R^{n-1})}+(w',\partial_nw_1)_{L_2(\rnp)}=(w',\partial_nw_1)_{L_2(\rnp)},
\endaligned
$$
where we used Proposition 6.1 in a similar way, and at the end used
that $\gamma _0w_1=0$, cf.\ (6.8).
Finally, taking the contributions from $\widehat P_2$ together, we get
$$
\multline
\int _{\rnp}\widehat P_2u\,{\partial_n \bar u'}\, dx+\int
_{\rnp}\partial_nu\,\overline
{\widehat P_2^* u'}\, dx
=(w_1,\partial_nw')_{L_2(\rnp)}+(\partial_nw_1,w')_{L_2(\rnp)}\\
=-(\gamma _0w_1,\gamma _0w')_{L_2(\R^{n-1})}=0,
\endmultline
$$
using again that $\gamma _0w_1=0$.

It is found in a similar way, using that $\widehat F_-^*$ is of
plus-type, that $\widehat P_3$ contributes with zero. \qed

\enddemo

\subsubhead 6.2 An integration by parts formula for $P$ \endsubsubhead

To extend the formula to the original operator $P$, we shall show that
 $P'=P-\widehat P$ (cf.\ Theorem 5.8 $1^\circ$) gives a zero boundary contribution.

\proclaim{Lemma 6.3} Let $a>0$ and let $S=\Op(s(\xi ))$, where $s(\xi )$ is
$O(\ang\xi ^{2a-1})$. Then
$$
\int_{\rnp}Su\,\partial_n\bar
u'\,dx+\int_{\rnp}\partial_nu\,\overline{S^*u'}\,dx =0,\tag6.10
$$
for any $u,u'\in \dot H^a(\crnp)$.
\endproclaim

  \demo{Proof} Since $u\in \dot H^a(\crnp)$, $r^+Su\in \ol
  H^{1-a}(\rnp)$; moreover $\partial_nu'\in \dot
  H^{a-1}(\crnp)$, so we can write the first integral  as
$$
\ang{r^+Su,\partial_nu'}_{\ol H^{1-a}(\rnp), \dot H^{a-1}(\crnp)}.
$$
Approximate $u'$ in $\dot H^a(\crnp)$ by a sequence of functions
$\varphi _k\in C_0^\infty (\rnp)$, $k\in\N$; then 
$$
\ang{r^+Su,\partial_n\varphi _k}_{\ol H^{1-a}, \dot
H^{a-1}}=-\ang{r^+\partial_nSu,\varphi _k }_{\ol H^{-a}, \dot
H^{a}}\to
-\ang{r^+\partial_nSu,u'}_{\ol H^{-a}, \dot H^{a}}. 
$$
With a similar argument for the second integral, we have
$$
\aligned
\ang{r^+Su,\partial_nu}_{\ol H^{1-a}, \dot
H^{a-1}}&+\ang{\partial_nu,r^+S^*u'}_{\dot H^{a-1}, \ol
H^{1-a}}\\
&=-\ang{r^+\partial_nSu,u' }_{\ol H^{-a}, \dot
H^{a}}-\ang{u,r^+\partial_nS^*u}_{\dot H^{a},\ol H^{-a}}\\
&=-\ang{r^+\partial_nSu,u' }_{\ol H^{-a}, \dot
H^{a}}+\ang{u,r^+(\partial_nS)^*u}_{\dot H^{a},\ol H^{-a}}=0,
\endaligned
$$
since $\partial_nS^*=S^*\partial_n=-(\partial_nS)^*$, and it is
well-known that the operator $S_1=\partial_nS$ of order $2a$ satisfies
$\ang{r^+S_1u,u'}=\ang{u,r^+S_1^*u'}$ for $u,u'\in \dot H^a(\crnp)$.
 \qed
\enddemo

We can then conclude:

\proclaim{Theorem 6.4} Let $P$, $\mu $, $\mu '$ be as in Theorem {\rm 6.2}.
  For $u\in \E_\mu (\crnp)\cap \E'(\rn)$, $u'\in \E_{\bar\mu '} (\crnp)\cap \E'(\rn)$,
 there holds 
$$
\aligned 
\int_{\rnp}  Pu\,\partial_n\bar
u'\,dx&+\int_{\rnp}\partial_nu\,\overline{  P^*u'}\,dx\\
&=\Gamma (\mu +1){\Gamma(\mu '+1)}\int_{{\Bbb
R}^{n-1}}s_0\gamma _0(u/x_n^{\mu })\,{\gamma _0(\bar u'/x_n^{\mu '})}\,
dx',
\endaligned
\tag6.11
$$
where $s_0=e^{-i\pi \delta }p(0,1)$. The formula extends to $u\in H^{\mu (t)}(\crnp)$,  $u'\in H^{\bar\mu
'(t')}(\crnp)$, for  $t>\operatorname{Re}\mu +\frac12$,
$t'>\operatorname{Re}\mu '+\frac12$, with $t,t'\ge a$.

The integrals over $\rnp$ are interpreted as dualities (as in
Proposition {\rm 6.1} and Lemma {\rm 6.3}) when $\operatorname{Re}\mu $ or
$\operatorname{Re}\mu '\le 0$, and when extended to general $u,u'$.  
\endproclaim

\demo{Proof} Recall that $P=\widehat P+P'$, where $P'=\Op(p')$,
  $|p'(\xi )|\le C\ang{\xi }^{2a-1}$ (cf.\ Theorem 5.8 $1^\circ$).
We
  have the identities (6.10) with $S=P'$ and (6.4) for $u\in H^{\mu (t)}(\crnp)$, $u'\in H^{\bar\mu
'(t')}(\crnp)$ with $t,t'\ge a$, $t>\operatorname{Re}\mu +\frac12$,
$t'>\operatorname{Re}\mu '+\frac12$.
Adding the identities for $\widehat P$ and $P'$ we obtain (6.11). It holds a fortiori for $u\in
\E_\mu (\crnp)\cap \E'(\rn)$, $u'\in \E_{\bar\mu '} (\crnp)\cap \E'(\rn)$. \qed

\enddemo

\example{Example 6.5}
The theorem applies in particular to $L=\Op(\Cal L(\xi ))$ studied in
(3.5)--(3.6) and
Example 5.9, showing that
$$
\aligned 
\int_{\rnp}  Lu\,\partial_n\bar
u'\,dx&+\int_{\rnp}\partial_nu\,\overline{  L^*u'}\,dx\\
&=\Gamma (\mu +1){\Gamma(\mu '+1)}\int_{{\Bbb
R}^{n-1}}|\Cal L(0,1)|\gamma _0(u/x_n^{\mu })\,{\gamma _0(\bar u'/x_n^{\mu '})}\,
dx',
\endaligned
\tag6.12
$$
The value $s_0=|\Cal L(0,1)|=(\Cal A(0,1)^2+\Cal B(0,1)^2)^{\frac12}$ is
found in (3.17).

This result was proved in \cite{DRSV21, Prop.\ 1.4} by completely different, real
methods, for $\mu \in \,]0,2a[\,\cap \,]2a-1,1[\,$.

The result is one of the key ingredients in the proof of integration
by parts formulas for operators $L$
on bounded domains $\Omega \subset\rn$ in \cite{DRSV21}, where $\mu
(\nu )$ varies as the normal $\nu $ varies along the boundary. It
would be interesting to extend this knowledge to general strongly
elliptic operators $P$ on bounded domains by similar applications of Theorem 6.5.
\endexample
\example{Example 6.6} Here is an example of an application to a
nonelliptic operator satisfying Assumption 3.2. Let 
$$
P=|\partial_1+\partial_2|^m+|\partial_3|^{m-1}\partial_3,\text{ with
symbol }p(\xi )=|\xi _1+\xi _2|^m+i\operatorname{sign}\xi _3|\xi _3|^m
$$
on ${\Bbb R}^3$, for some  $1<m<2$. For $\R^3_+=\{x_3> 0\}$ we have the normal
$\nu =(0,0,1)$, where
$$
p(0,0,1)=i,\; p(0,0,-1)=-i,\text{ so (3.1) holds with }m-2\mu =1,
$$
i.e., $\mu =(m-1)/2$, $\mu '=(m+1)/2$. Then by Theorem 6.4,
$$
\int_{x_3>0}(Pu\, \bar v-u\,\overline{P^*v})\, dx=\Gamma
(\tfrac{m}2+\tfrac12)\Gamma (\tfrac{m}2+\tfrac32)\int_{x_3=0}\gamma _0\bigl(\tfrac
u{x_3^{(m-1)/2}}\bigr) \gamma _0\bigl(\tfrac
{\bar v}{x_3^{(m+1)/2}}\bigr)\,dx',
$$
for functions $u\in x_3^{(m-1)/2}\Cal S(\overline{\R}^3_+)$, $v\in x_3^{(m+1)/2}\Cal S(\overline{\R}^3_+)$.

The halfspace $\{x_2>0\}$ has the normal $\nu
'=(0,1,0)$ and 
$$
p(\nu ')=1,\; p(-\nu ')=1,\text{ so (3.1) holds with }m-2\mu =0,
$$
i.e., $\mu =m/2$, $\mu '=m/2$. Here by Theorem 6.4,
$$
\int_{x_2>0}(Pu\, \bar v-u\,\overline{P^*v})\, dx=\Gamma
(\tfrac{m}2+1)^2\int_{x_2=0}\gamma _0\bigl(\tfrac
u{x_2^{m/2}}\bigr) \gamma _0\bigl(\tfrac
{\bar v}{x_2^{m/2}}\bigr)\, dx_1dx_3,
$$
for functions with a factor $ x_2^{m/2}$. 
\endexample

\subhead 7. Large solutions and a halfways Green's formula\endsubhead

\subsubhead 7.1 Large solutions, a nonhomogeneous Dirichlet problem \endsubsubhead

Let $P$ satisfy Assumption 3.1, and assume $\operatorname{Re}\mu >0$.
Along with the homogeneous Dirichlet problem (5.25), one can consider a
nonhomogeneous local Dirichlet problem if  the scope is expanded to allow
so-called "large solutions", behaving like $x_n^{\mu -1}$ near the
boundary of $\rnp$; such solutions blow up at the boundary when $\operatorname{Re}\mu <1$.
Namely, one can pose the nonhomogeneous Dirichlet problem 
$$
r^+Pu=f \text{ on }\rnp ,\quad \gamma _0(u/x_n^{\mu -1} )=\varphi  
\text{ on }\R^{n-1} ,\quad \supp u\subset \crnp.\tag7.1
$$

Problem (7.1) was studied earlier for operators satisfying the
$a$-transmission property in \cite{G15,G14} (including the fractional
Laplacian $(-\Delta )^a$), and a halfways Green's formula was shown in
\cite{G18, Cor.\ 4.5}. The problem (7.1) for the fractional Laplacian, and the
halfways Green's formula --- with applications to solution formulas
--- were also
studied in
Abatangelo \cite{A15} (independently of
\cite{G15});  the boundary condition there is given
in a less explicit way except when $\Omega $ is a ball. 
There have been further studies of such problems, see e.g.\ Abatangelo,
Gomez-Castro and Vazquez \cite{AGV21} and its references.

Note that the boundary condition in (7.1) is {\it local}. There is a different
problem which is also regarded as a nonhomogeneous Dirichlet problem, 
namely to
prescribe nonzero values of $u$ in the exterior of $\Omega $; it has
somewhat different solution spaces (a link between this and the
homogeneous Dirichlet problem is described in \cite{G14}).

For the general  operators $P$ considered here, we shall now show that problem (7.1) has a good sense for $u\in H^{(\mu
-1)(t)}(\crnp)$ with suitable $t$.

More precisely, since $P$ also satisfies the principal $(\mu
-1)$-transmission condition (as remarked after Definition 2.1),  Theorem 5.8 $1^\circ$
can be applied with $\mu $ replaced by $\mu -1$, implying that $r^+P$
maps
$$
r^+ P\colon  H^{(\mu -1)(t) }(\crnp)\to \ol H^{t -2a }(\rnp)\text{ for
}\operatorname{Re}\mu -\tfrac 32<t<\operatorname{Re}\mu
+\tfrac12. \tag7.2
$$
This is also valid in the case where $P$ is only assumed to satisfy
Assumption 3.2.

From Theorem 5.4 we have (note that  $\operatorname{Re}\mu -1>-1$)
$$
H^{(\mu -1)(t) }(\crnp)\cases =\dot H^t (\crnp)\text{ when
}-\frac32<t-\operatorname{Re}\mu <-\frac12
,\\
\subset \dot H^t(\crnp)+x_n^{\mu -1}e^+\ol H^{t-\operatorname{Re}\mu +1} (\rnp)\text{ when
}-\frac12<t-\operatorname{Re}\mu <\frac12.
\endcases \tag7.3
$$
When $t-\operatorname{Re}\mu >-\frac12$, the weighted
boundary value is well-defined, cf.\ (5.20):
$$
\gamma _0^{\mu -1}u\equiv\gamma _0(\Xi _+^{\mu -1}u)=\Gamma (\mu )\gamma _0(u/x_n^{\mu -1})\in H^{t-\operatorname{Re}\mu +\frac12}(\R^{n-1}).\tag7.4
$$

The following regularity result holds for the nonhomogeneous Dirichlet
problem:

\proclaim{Theorem 7.1} Let $P$ satisfy Assumption {\rm 3.1} with
$\operatorname{Re}\mu >0$, and let $-\frac12<t-\operatorname{Re}\mu
<\frac12$. When $f\in \ol H^{t-2a}(\rnp)$ and $\varphi \in
H^{t-\operatorname{Re}\mu +\frac12}(\R^{n-1})$ are given, and $u$
solves the nonhomogeneous Dirichlet problem {\rm (7.1)}
with $u\in H^{(\mu -1)(\sigma ) }(\crnp)$ for some  $-\frac12<\sigma -\operatorname{Re}\mu
<\frac12$, then in fact $u\in H^{(\mu -1)(t ) }(\crnp)$.
\endproclaim

\demo{Proof} It is known from \cite{G15, Th.\ 7.1} that $H^{\mu (\sigma )
}(\crnp)$ is a closed subspace of $ H^{(\mu
-1)(\sigma ) }(\crnp)$,  equal to the set of $v\in H^{(\mu
-1)(\sigma ) }(\crnp)$ for which $\gamma _0(v/x_n^{\mu -1})=0$.
From the given $\varphi $ we define
$$
w=\Gamma (\mu )\Xi _+^{-\mu +1}e^+K_0\varphi \in H^{(\mu -1)(t ) }(\crnp),
$$
where $K_0$ is the standard Poisson operator $\varphi \mapsto
K_0\varphi =\Cal F^{-1}_{\xi '\to x'}[\hat \varphi (\xi ')e^{-\ang{\xi
'}x_n}]$, $x_n>0$. Then in view of (7.4),
$$
\gamma _0(w/x_n^{\mu -1})=\Gamma (\mu )^{-1}\gamma _0(\Xi _+^{\mu -1}w)=
       \gamma _0(\Xi _+^{\mu -1}\Xi _+^{1-\mu }e^+K_0\varphi )
=\gamma _0K_0\varphi =\varphi ,
$$ so that $v=u-w$ solves (7.1) with $f$ replaced by
$f-r^+Pw\in \ol H^{t-2a}(\rnp)$, $\varphi $ replaced by 0. This is
a homogeneous Dirichlet problem as in (5.25). Since   $v\in H^{(\mu
-1)(\sigma ) }(\crnp)$ with $\gamma _0(v/x_n^{\mu -1})=0$, it is in $
H^{\mu (\sigma ) }(\crnp)$. It then follows from Theorem 5.8 that
$v\in H^{\mu (t ) }(\crnp)$, and hence $u=v+w\in H^{(\mu -1)(t )
}(\crnp)$. \qed 
\enddemo

For the hatted version $\widehat P$ there is even an existence and
uniqueness result in these spaces.

\proclaim{Theorem 7.2} Let $P$ satisfy Assumption {\rm 3.1} with
$\operatorname{Re}\mu >0$, and let $-\frac12<t-\operatorname{Re}\mu
<\frac12$. Then $r^+\widehat P$ together with $\gamma _0^{\mu -1}$ defines a homeomorphism:
$$
\{r^+\widehat P,\gamma _0^{\mu -1}\}\colon  H^{(\mu -1)(t ) }(\crnp)\simto \ol H^{t-2a}(\rnp)\times
H^{t-\operatorname{Re}\mu +\frac12}(\R^{n-1}).\tag7.5
$$
\endproclaim

\demo{Proof} The forward mapping properties are accounted for
above. The existence of a unique solution $u\in H^{(\mu -1)(t )
}(\crnp) $ of 
$$
r^+\widehat Pu=f\text{ on }\rnp,\quad \gamma _0^{\mu -1}u=\varphi
\text{ on }\R^{n-1},\quad \supp u\subset \crnp,\tag7.6
$$
 for given $f\in \ol H^{t-2a}(\rnp)$, $\varphi \in
H^{t-\operatorname{Re}\mu +\frac12}(\R^{n-1})$, is shown as in Theorem
7.1, now referring to Theorem 5.5 instead of Theorem 5.8. 
\qed 
\enddemo

These theorems  show that  $H^{(\mu -1)(t ) }(\crnp)$ is the correct domain space for
the nonhomogeneous Dirichlet problem, at least in the small range 
$-\frac12<t-\operatorname{Re}\mu
<\frac12$. Recall that  $\E_{\mu -1}(\crnp)\cap \Cal E'(\rn)$ and $e^+x_n^{\mu -1}\Cal S(\crnp)$ are
dense subsets of $H^{(\mu -1)
(t)}(\crnp)$ for all $t>\operatorname{Re}\mu -\frac32$.

\subsubhead 7.2 An integration by parts formula involving the nontrivial Dirihlet trace \endsubsubhead

We now show a ``halfways Green's formula'', where one factor $u$
is in the domain of the nonhomogeneous Dirichlet problem for $P$ and the other
factor $v$ is in the domain of the homogeneous Dirichlet
problem for $P^*$:

\proclaim{Theorem 7.3} Let  $P$ satisfy Assumption {\rm 3.2},
 and assume moreover that $0<\operatorname{Re}\mu <a+\frac12$, $ \operatorname{Re}\mu ' >0$.
 
 For $u\in \E_{\mu -1}(\crnp)\cap \E'(\rn)$ and  $v\in \E_{\bar\mu '} (\crnp)\cap \E'(\rn)$,
 there holds 
$$
\int_{\rnp} Pu\,\bar v\,dx-\int_{\rnp}u\,\overline{P^*v}\,dx=-\Gamma (\mu ){\Gamma(\mu '+1)}\int_{{\Bbb
R}^{n-1}}s_0\gamma _0(u/x_n^{\mu -1 })\,{\gamma _0(\bar v/x_n^{\mu '})}\,
dx',
\tag7.7
$$
where $s_0=e^{-i\pi \delta }p(0,1)$.
The formula extends to $u\in H^{(\mu -1)(t)}(\crnp)$ with
$t>\operatorname{Re}\mu -\frac12$, $v\in H^{\bar\mu '(t')}(\crnp)$
with $t'>\operatorname{Re}\mu '+\frac12$.

The left-hand side is interpreted as follows, for small $\varepsilon >0$:
$$
\ang{r^+Pu,v}_{\ol H^{-\operatorname{Re}\mu '-\frac12+\varepsilon}(\rnp),\dot H^{\operatorname{Re}\mu '+\frac12-\varepsilon}(\crnp)}-\ang{u,P^*v}_{\dot H^{\operatorname{Re}\mu -\frac12-\varepsilon}(\crnp),\ol H^{- \operatorname{Re}\mu
+\frac12+\varepsilon}(\rnp)}.\tag 7.8
$$

\endproclaim

\demo{Proof} We shall show how the result can be derived from Theorem
6.4. Let $u\in \E_{\mu -1}(\crnp)\cap \E'(\rn)$ and  $v\in \E_{\bar\mu '}
(\crnp)\cap \E'(\rn)$. As shown in \cite{G15, p.\ 494}, there exist
functions $U$ and $u_1$ in $\E_{\mu }(\crnp)\cap \E'(\rn)$ such that
$u=\partial_nU+u_1$.

In terms of the Hilbert spaces: When $u\in H^{(\mu -1)(t)}(\crnp)$ with
$|t-\operatorname{Re}\mu |<\frac12$, let  $z=r^+\Xi _+^{\mu -1}u\in
\ol H^{t-\operatorname{Re}\mu +1}(\rnp)$, then (denoting $\Op(\ang{\xi '})=\ang{D'}$)
$$
\aligned
u&=\Xi _+^{-\mu +1}e^+z=(\ang{D'}+\partial_n)\Xi _+^{-\mu
}e^+z
=u_1+\partial_nU,\text{ with}\\
u_1&=\ang{D'}\Xi _+^{-\mu
}e^+z\in \ang{D'}H^{\mu (t+1)}(\crnp)\subset H^{\mu (t)}(\crnp),\\ 
U&=\Xi _+^{-\mu }e^+z \in H^{\mu (t+1)}(\crnp),\quad \partial_nU\in H^{\mu (t)}(\crnp). 
\endaligned \tag 7.9
$$
Here $ H^{\mu (t)}(\crnp)=\dot H^t(\crnp)$ since  $|t-\operatorname{Re}\mu
|<\frac12$. Moreover, when $t=\operatorname{Re}\mu
-\frac12+\varepsilon $ for a small $\varepsilon >0$, then 
$$
r^+Pu=r^+Pu_1+r^+P\partial_nU=r^+Pu_1+\partial_nr^+PU
$$
where both terms are in $\ol H^{t-2a}(\rnp)=\ol
H^{\operatorname{Re}\mu-\frac12 +\varepsilon -2a}(\rnp)=\ol
H^{\,-\operatorname{Re}\mu '-\frac12+\varepsilon }(\rnp)$; we here use
Theorem 5.8 $1^\circ$. 

For $v$, we note that when $v\in H^{\bar\mu '(t')}(\crnp)$ with
$t'=\operatorname{Re}\mu '+\frac12+\varepsilon $, then
$$
\aligned
v&\in H^{\bar\mu '(t')}(\crnp)=\Xi _+^{-\bar\mu '}e^+\ol
H^{t'-\operatorname{Re}\mu '}(\rnp)\subset \Xi _+^{-\bar\mu '}
\dot H^{\frac12-\varepsilon }(\crnp)= \dot H^{\operatorname{Re}\mu
'+\frac12-\varepsilon }(\crnp),\\
r^+P^*v&\in \ol H^{t'-2a}(\rnp)= \ol H^{\operatorname{Re}\mu '+\frac12+\varepsilon -2a}(\rnp)= \ol H^{\,-\operatorname{Re}\mu +\frac12+\varepsilon }(\rnp).
\endaligned
$$
Then the dualities in (7.8) are well-defined and serve as an
interpretation of the left-hand side in (7.7).

The formula (7.7) will first be proved for $u\in \E_{\mu -1}(\crnp)\cap \E'(\rn)$ and  $v\in \E_{\bar\mu '}
(\crnp)\cap \E'(\rn)$, and afterwards be
 extended by continuity to general $u,v$.
We use the decomposition (7.9), that leads to elements of
$\E_\mu (\crnp)$ for $t\to \infty $. When $u$ is supported in a ball
$\{|\xi |\le R\}$,
we can cut $u_1$ and $U$ down to have support in $\{|\xi |\le 2R\}$. 

Consider the contribution from $u_1$. Here there holds 
$$
\ang{r^+Pu_1,v}_{\ol H^{-a}, \dot H^a}-\ang{u_1,r^+P^*v}_{\dot H^a,
\ol H^{-a}}=0,\tag7.10
$$
when $u_1$ and $v$ are in $\dot H^a(\crnp)$, since $P$ is of order $2a$. This gives the
contribution 0 to (7.7) since $t=a$ is allowed in the definition of
$u_1$ (recall that $a>\operatorname{Re}\mu -\frac12$ by
hypothesis), and $t'\ge a$ holds for the values of $t'$ allowed in the definition of  $v$ (where $t'>\operatorname{Re}\mu '+\frac12=2a-\operatorname{Re}\mu +\frac12>2a-a=a$).
Thus $u_1$ contributes to the boundary integral with 0.

For the contribution from $\partial_nU$, we note that, writing
$U=x_n^\mu w$ for $x_n>0$, $w\in C^\infty (\crnp)$,
$$
\partial_nU=\partial_n(x_n^\mu w)=\mu x_n^{\mu
-1}w+x_n^\mu \partial_nw 
\text{ for }x_n>0,
$$
so the weighted boundary value for $x_n\to 0+$ satisfies (since
$x_n^\mu\partial_nw/x_n^{\mu -1}=x_n\partial_nw \to 0$)
$$
\gamma _0(\partial_nU/x_n^{\mu -1})=\mu \gamma _0w=\mu \gamma _0(U/x_n^{\mu }).\tag7.11
$$
Moreover, by a simple integration by parts,
$$
\ang{r^+P\partial_nU,v}=\ang{r^+\partial_nPU,v}=-\ang{r^+PU,\partial_nv},
$$
since $\gamma _0v=0$ because of $\operatorname{Re}\mu '>0$.
Thus, by use of Theorem 6.4 and (7.11),
$$
\aligned
\ang{r^+P\partial_nU,v}&-\ang{\partial_nU,r^+P^*v}=-\ang{r^+PU,\partial_nv}-\ang{\partial_nU,r^+P^*v}\\
&=-\Gamma
(\mu +1)\Gamma (\mu '+1)s_0\int_{\R^{n-1}}\gamma _0(U/x_n^{\mu })\gamma
_0(\bar v/x_n^{\mu '})\,dx'\\
&=-\Gamma
(\mu )\Gamma (\mu '+1)s_0\int_{\R^{n-1}}\gamma _0(\partial_nU/x_n^{\mu -1})\gamma
_0(\bar v/x_n^{\mu '})\, dx'.
\endaligned
$$
Since $u_1\in \E_\mu (\crnp)$, $\gamma _0(u_1/x_n^{\mu -1})=0$, so
$u_1$ can be added to $\partial_nU$ in the last integral. Adding also
(7.10) to the left-hand side, we find  (7.7).

Since the expressions depend continuously on $u,v$ in the presented
norms, the formula extends to the indicated spaces. \qed

\enddemo

\example{Example 7.4} Theorems 7.1 and 7.2 apply in particular to the
operator $L$ considered in (3.5)--(3.6) and Examples 5.9 and 6.5, when $\mu >0$ (this
holds automatically if $a\ge \frac12$, since $|\delta |<\frac12$). Theorem 7.3 applies to $L$
when $\mu $ and $\mu '>0$ (again automatically satisfied when $a\ge \frac12$).
\endexample

\example{Remark 7.5} The transmission spaces can also be defined in
terms of other scales of function spaces. The case of Bessel-potential
spaces $H^s_p$, $1<p<\infty $, is a main subject in our preceding
papers. There is also
the H\"older-Zygmund scale $C_*^s(\rn)$, coinciding with the
H\"older scale $C^s(\rn)$ when $s\in \rp\setminus \N$, with spaces
over $\rnp$ defined as in (5.1). Here since $C_*^{s+\varepsilon
}(\rn)\subset H^s(\rn)$, also $C_*^{\mu (s+\varepsilon
)}(\crnp)\subset H^{\mu (s)}(\crnp)$ for $\varepsilon >0$. (More details
on such spaces in our earlier papers, e.g.\ in \cite{G19}.) So the
results dealing with forward mapping properties of $r^+P$ have useful
consequences involving these spaces as well. Namely, Theorem 5.8 $1^\circ$
implies that $r^+P$ maps
$$
r^+ P\colon  C_*^{\mu (t+\varepsilon ) }(\crnp)\to \ol H^{t -2a }(\rnp)\text{ for }\operatorname{Re}\mu -\tfrac12<t<
\operatorname{Re}\mu +\tfrac32,
$$
and the integration by parts formulas in Sections 5 and 6 hold for
functions   in $ C_*^{\mu (t)}$-type spaces, for the same $t$.

In the opposite direction, an inclusion of an $H^s$-space in a H\"older spaces loses $n/2$ in the
regularity parameter, so does not give very good results. For better
regularity results, it
would be interesting to generalize the above theory to
$H^s_p$-spaces with general  $1<p<\infty $, possibly under further
hypotheses; this remains to be done. More smoothness than $C^1$ is needed for
a symbol $q(\xi )$ to be a Fourier multiplier in $L_p$ (some well-known
conditions are recalled in \cite{GK93, Sect.\ 1.3}). 
There is an extension of Vishik and Eskin's
work to $L_p$-based spaces by Shargorodsky \cite{S94}, which should be
useful. It is there
pointed out that \cite{E81, Lemma 17.1} shows how smoothness properties
carry over to the factors in the Wiener-Hopf factorization. 
\endexample

\enddocument